\documentclass[a4paper]{amsart}
\usepackage[top=35mm,bottom=15mm,left=20mm,right=20mm]{geometry}
\usepackage{graphicx} 
\usepackage{amsmath,amssymb,amsthm,mathtools}
\usepackage[colorlinks,citecolor=blue,urlcolor=blue,bookmarks=true]{hyperref}
\usepackage{comment}
\usepackage[shortlabels]{enumitem}
\usepackage[dvipsnames]{xcolor}

\newtheorem{theorem}{Theorem}
\newtheorem{proposition}{Proposition}
\newtheorem{lemma}{Lemma}
\newtheorem*{rem}{Remark}
\theoremstyle{definition}
\newtheorem*{dfn}{Definition}
\newtheorem{cor}{Corollary}

\title{Squared Bessel particle systems with weak repulsion}
\author[J.~Lenczewska]{Julia Lenczewska}
\address{Wroc{\l}aw University of Science and Technology,
Faculty of Pure and Applied Mathematics\\
	Wyb. {Wyspia\'{n}skiego} 27,
	50-370 {Wroc{\l}aw}, Poland}
    \email{julia.lenczewska@pwr.edu.pl}
\author[J.~Małecki]{Jacek Małecki}
\address{Wroc{\l}aw University of Science and Technology,
Faculty of Pure and Applied Mathematics\\
	Wyb. {Wyspia\'{n}skiego} 27,
	50-370 {Wroc{\l}aw}, Poland}
\email{jacek.malecki@pwr.edu.pl}

\author[S.~Sołtysiak]{Szymon Sołtysiak}
\address{Wroc{\l}aw University of Science and Technology,
Faculty of Pure and Applied Mathematics\\
	Wyb. {Wyspia\'{n}skiego} 27,
	50-370 {Wroc{\l}aw}, Poland}
\email{szymon.soltysiak@pwr.edu.pl}
\date{\today}

\subjclass[2020]{Primary 60H10; Secondary 60J60, 60K35}

\keywords{Squared Bessel particle systems,
singular stochastic differential equations,
weak repulsion,
pathwise uniqueness,
non-sticky solutions,
particle collisions,
hitting times}

\begin{document}

\begin{abstract}
    We study squared Bessel particle systems with repulsion parameter $\kappa>0$ and dimension parameter $\delta\in\mathbb R$. We prove that the singular equation has a unique global strong solution for every initial state in the closed Weyl chamber, including states with coinciding particles. The solutions have locally integrable interaction terms and spend zero Lebesgue time at pairwise collisions. We also consider the equation in which the interaction terms are set to zero at collisions. At certain integer values of $\delta$, this equation admits infinitely many solution laws from suitable initial configurations. In these solutions, a block of particles can remain at zero for any prescribed time before being released. For positive initial states, we determine when zero is reached and when particles enter the negative half-line. From distinct initial positions, collisions occur almost surely for $0<\kappa<1$ and do not occur for $\kappa\ge1$. Collisions of three or more particles are excluded away from zero. We give an example of a triple collision at zero. We obtain two-sided bounds for the tails of the first-zero and first-collision times. The bounds give matching powers for the first-zero time when it is finite, and for the first-collision time in nonnegative systems when $0<\kappa<1$.
\end{abstract}

\maketitle

\section{Introduction}

We consider squared Bessel particle systems with an arbitrary positive repulsion parameter. Let $N\ge2$, $\kappa>0$ and $\delta\in\mathbb R$. The particles are ordered and take values in the closed Weyl chamber
    \begin{equation*}
        \mathcal C_A=\{x\in\mathbb R^N:x_1\le\cdots\le x_N\}.
    \end{equation*}
We study the system
    \begin{equation}
        \label{eq:A-particle-genuine}
        dX_i(t) = 2\sqrt{|X_i(t)|}\,dB_i(t)+\kappa\left(\delta +\sum_{j\ne i} \frac{|X_i(t)|+|X_j(t)|}{X_i(t)-X_j(t)}\right)dt, \qquad i=1,\ldots,N,
    \end{equation}
where $B_1,\ldots,B_N$ are independent standard Brownian motions. We allow every initial point $x\in\mathcal C_A$, also when some coordinates are equal or equal to zero. Our aim is to construct global solutions and to describe their collisions and their behaviour at zero. For small values of $\kappa$, collisions may occur after the initial time. At a collision at zero, both the numerator of the interaction term and the diffusion coefficient vanish. This case requires a separate analysis. We also use the equation
    \begin{equation}
        \label{eq:A-particle-indicator}
        dX_i(t) = 2\sqrt{|X_i(t)|}\,dB_i(t)+\kappa\left(\delta +\sum_{j\ne i} \frac{|X_i(t)|+|X_j(t)|}{X_i(t)-X_j(t)} \mathbf1_{\{X_i(t)\ne X_j(t)\}}\right)dt, \qquad i=1,\ldots,N.
    \end{equation}
Here each interaction term is defined to be zero when the two coordinates are equal. We call \eqref{eq:A-particle-indicator} the \emph{indicator equation}. This convention can give solutions which spend positive time at collisions. For \eqref{eq:A-particle-genuine} we require the quotients to be defined almost everywhere in time and locally absolutely integrable. Thus the solution we seek spends zero Lebesgue time at every pairwise collision. The precise definitions are given below.

\bigskip

The study of interacting eigenvalues goes back to Dyson \cite{Dyson1962}. For squared Bessel particles, the matrix model is the Wishart process introduced by Bru \cite{Bru1991}. In the usual real Wishart normalization, the eigenvalue equation agrees with \eqref{eq:A-particle-genuine} on the nonnegative chamber for $\kappa=1$. The complex Laguerre case corresponds to $\kappa=2$. K\"onig and O'Connell \cite{KonigOConnell2001} represented the latter eigenvalue process as a Doob transform of independent squared Bessel processes killed at their first collision. This is related to the determinantal formula of Karlin and McGregor \cite{KarlinMcGregor1959} and to the study of noncolliding processes in Weyl chambers by Grabiner \cite{Grabiner1999}. Matrix models for general beta ensembles were constructed by Dumitriu and Edelman \cite{DumitriuEdelman2002}. These results explain the random matrix origin of the interaction, but the signed equation with arbitrary parameters is a separate stochastic differential equation problem. In one dimension, the use of $\sqrt{|x|}$ allows squared Bessel equations with negative dimensions and signed initial values; see G\"oing-Jaeschke and Yor \cite{GoingJaeschkeYor2003}.

Strong existence for singular particle equations was studied by Graczyk and Ma\l ecki in \cite{GraczykMalecki2013,GraczykMalecki2014}. The first paper gives a multidimensional Yamada-Watanabe argument with applications to eigenvalue equations. The second constructs noncolliding systems from arbitrary ordered initial points, including colliding ones, by using elementary symmetric polynomials. Its applications include nonnegative squared Bessel systems for $\kappa\ge1$ and $\delta\ge N-1$. Thus coinciding initial coordinates are already allowed in the earlier noncolliding theory.

The closest earlier paper on the signed model is Graczyk and Ma\l ecki \cite{GraczykMalecki2019}. It treats \eqref{eq:A-particle-indicator} for $\kappa=1$, every real $\delta$ and every ordered initial point. The authors construct the unique noncolliding strong solution and describe its behaviour at zero and its separation into positive and negative groups. They also study uniqueness for the full indicator equation, including exceptional parameters for which other solutions exist, and characterize nonnegative solutions. In particular, pinned zero coordinates and nonuniqueness of the indicator equation already appear there. The present paper keeps the signed state space and allows every $\kappa>0$, including the range in which the selected solution itself can have collisions.

Colliding systems with positive repulsion have a longer history. C\'epa and L\'epingle \cite{CepaLepingle1997} constructed Brownian particle systems with electrostatic repulsion, also when the repulsion is too weak to prevent pair collisions. They later proved that collisions of three or more particles do not occur in that model \cite{CepaLepingle2007}. Radial Dunkl processes give a related description in terms of root systems. Demni \cite{Demni2009} proved strong existence and uniqueness for strictly positive multiplicities and studied the hitting of chamber walls. These results allow collisions, but their direct application to squared Bessel particles does not cover all values of $\delta$ or the passage to signed coordinates.

For nonnegative beta-Wishart systems, weak repulsion was studied in detail by Jourdain and Kahn \cite{JourdainKahn2022}. Their parameters correspond to ours by $\beta=\kappa$, $\alpha=\kappa\delta$ and $\gamma=0$. Write
    \begin{equation*}
        d_A=\kappa(\delta-N+1).
    \end{equation*}
For $0<\kappa<1$, $d_A>0$ and initial states with at most one zero coordinate, their construction gives a pathwise unique strong solution up to the first time that two particles are simultaneously zero. It is global when $d_A\ge1-\kappa$. When $0<d_A<1-\kappa$, that double-zero time is finite. They also study pair collisions and the approach of groups of particles to zero. Continuation past the double-zero time is supplied by the later nonnegative construction in Ma\l ecki \cite[Proposition~18]{Malecki2026}, which gives a global non-sticky strong solution for every $\kappa>0$ and $\delta>N-1$.

The general theory in \cite{Malecki2026} also applies to the signed indicator equation considered here. Proposition~19 gives global weak solutions of the indicator equation and proves pathwise uniqueness among non-sticky solutions. The remaining question is the construction of such a solution at degenerate zero configurations; this issue is discussed in Remark~17 of that paper. We resolve it by analysing collision occupation and constructing continuations from zero blocks. Recent results of Andraus, Hufnagel and Ma\l ecki \cite{AndrausHufnagelMalecki2025} describe collision types and times for more general particle systems with strictly positive diffusion and repulsion coefficients. The coefficients in \eqref{eq:A-particle-genuine} do not satisfy that positivity assumption at zero.

There are also explicit results on hitting times in the matrix and Dunkl settings. Donati-Martin, Doumerc, Matsumoto and Yor \cite{MR2105711} obtained formulas for the first zero of the determinant of a real Wishart process. Demni \cite{Demni2007} obtained the corresponding formulas for complex Laguerre processes. The first collision of independent squared Bessel processes, before the Doob transform, was studied in \cite{KonigOConnell2001}. For radial Dunkl processes, Demni \cite{MR2470522} computed distributions of the first hitting time of the boundary of a Weyl chamber. These are direct precedents for the hitting time questions studied below.

The distinction between the stopping times is relevant here. Starting from $0<x_1<\cdots<x_N$, the exit time from the positive ordered chamber is the minimum of the first zero time and the first pair collision time. Before that exit, $Y_i=\sqrt{X_i}$ satisfies a type-$B_N$ radial Dunkl equation with multiplicities
    \begin{equation*}
        k_0=\frac{d_A-1}{2},\qquad k_1=\frac{\kappa}{2}.
    \end{equation*}
The usual nonnegative multiplicity range requires $d_A\ge1$. Our first-zero time is defined through positive collisions, and our first-collision time is defined through individual visits to zero and changes of sign. Thus neither of these times is in general the positive-chamber exit time.

\subsection{Solutions and main results}

A continuous $\mathcal C_A$-valued process $X$ is called \emph{non-sticky} if, almost surely, for every $T>0$ and every $i<j$,
    \begin{equation*}
        \int_0^T\mathbf1_{\{X_i(t)=X_j(t)\}}\,dt=0.
    \end{equation*}
A solution of \eqref{eq:A-particle-genuine} consists of a continuous adapted $\mathcal C_A$-valued process $X$ and an $N$-dimensional Brownian motion $B$ on a filtered probability space, such that the equation holds in integral form. We require, almost surely, that each pairwise quotient be defined for Lebesgue-almost every time and that, for every $T>0$ and $i<j$,
    \begin{equation*}
        \int_0^T \frac{|X_i(t)|+|X_j(t)|}{|X_i(t)-X_j(t)|}\,dt<\infty.
    \end{equation*}
Values assigned on a Lebesgue-null set do not affect these integrals. For \eqref{eq:A-particle-indicator}, we impose the same continuity, adaptation and integral-equation requirements, with local absolute integrability of the interaction terms extended by zero at collisions. We do not require a solution of the indicator equation to be non-sticky.

Under this convention, every solution of \eqref{eq:A-particle-genuine} is non-sticky: its unextended quotients are undefined at every collision, also when both coordinates are zero. Conversely, every non-sticky solution of the indicator equation solves \eqref{eq:A-particle-genuine}. Non-stickiness permits collisions on sets of zero Lebesgue measure. It concerns pairs of coordinates and does not exclude a single coordinate which stays at zero.

Our first result gives the global solution under this convention.

\begin{theorem}[Global existence and uniqueness]
    \label{thm:A-genuine-existence-uniqueness}
    Let $N\ge2$, $\kappa>0$ and $\delta\in\mathbb R$. For every $x\in\mathcal C_A$, the system
        \begin{equation*}
        dX_i(t) = 2\sqrt{|X_i(t)|}\,dB_i(t)+\kappa\left(\delta +\sum_{j\ne i} \frac{|X_i(t)|+|X_j(t)|}{X_i(t)-X_j(t)}\right)dt, \qquad i=1,\ldots,N,
    \end{equation*}
    with $X(0)=x$, where $B_1,\ldots,B_N$ are independent standard Brownian motions, has a pathwise unique global strong solution taking values in $\mathcal C_A$.
\end{theorem}

For $\kappa=1$, this agrees with the noncolliding construction of \cite{GraczykMalecki2019}. In relation to \cite{Malecki2026}, the additional step is the construction of a global non-sticky continuation at every zero configuration, for arbitrary signed initial states and real $\delta$. Weak existence of the indicator equation and uniqueness in the non-sticky class are supplied by that general theory. The proof given here first uses a local-time argument to identify where positive collision occupation can occur. We then construct the required continuations from one-sided squared Bessel systems. The pathwise uniqueness argument uses Tanaka's formula and the monotonicity of the pair interaction.

The distinction between the two equations is necessary. Put
    \begin{equation*}
        \mathcal R_N=\mathbb Z\cap[2-N,N-2].
    \end{equation*}
If $\delta\notin\mathcal R_N$, every solution of the indicator equation is non-sticky. At suitable initial configurations with $\delta\in\mathcal R_N$, a block of zero coordinates can instead remain pinned for any prescribed finite time and then follow the non-sticky solution, or remain pinned forever. These choices give distinct strong solutions driven by the same Brownian motion and uncountably many solution laws. Such zero-block nonuniqueness is known for $\kappa=1$ from \cite{GraczykMalecki2019}; here the construction is given for every $\kappa>0$.

We next describe the first encounter with zero. For a strictly positive ordered initial point, let
    \begin{equation*}
        T_0=\inf\{t\ge0:X_1(t)=0\},\qquad T_-=\inf\{t\ge0:X_1(t)<0\}.
    \end{equation*}
As usual, the infimum of the empty set is $\infty$. The parameter which determines the behaviour at zero is $d_A=\kappa(\delta-N+1)$.

\begin{theorem}[Behaviour at zero]
    \label{thm:A-complete-zero-trichotomy}
    Let $N\ge2$, $\kappa>0$, $\delta\in\mathbb R$, and let $X$ be the solution of \eqref{eq:A-particle-genuine} started from $0<x_1<\cdots<x_N$. Then the following statements hold almost surely.
        \begin{enumerate}[(i)]
            \item If $d_A\ge2$, then $T_0=\infty$.
            \item If $0\le d_A<2$, then $T_0<\infty$ and $T_-=\infty$.
            \item If $d_A<0$, then $T_-=T_0<\infty$ and
                \begin{equation*}
                    \inf\{s>0:X_1(T_0+s)<0\}=0.
                \end{equation*}
        \end{enumerate}
\end{theorem}

The thresholds $0$ and $2$ are the familiar squared Bessel thresholds. The $\kappa=1$ signed behaviour is contained in \cite{GraczykMalecki2019}, and the nonnegative beta-Wishart boundary problem is studied in \cite{JourdainKahn2022}. The theorem gives the classification for every $\kappa>0$, also when collisions occur before $T_0$ or several coordinates reach zero together. The last case concerns immediate entrance after $T_0$; continuity still gives $X_1(T_0)=0$. One ingredient of the proof is the equation for $E=\prod_iX_i$ before $T_0$:
    \begin{equation*}
        dE=2\sqrt{EG}\,dW+d_A G\,dt, \qquad G=\sum_{i=1}^N\prod_{j\ne i}X_j.
    \end{equation*}
It identifies $E$ with a squared Bessel process after a time change. The analysis of a zero block is then used to determine what happens at and after the first zero.

The same product equation also gives quantitative bounds. For $a>0$, $t>0$ and $z\ge0$, write
    \begin{equation*}
        F_a(t,z)=\frac1{\Gamma(a)} \int_0^{z/(2t)}u^{a-1}e^{-u}\,du.
    \end{equation*}
For a fixed initial point, $f(t)\asymp g(t)$ means that their ratio is bounded above and below by positive constants for all sufficiently large $t$.

\begin{theorem}[First-zero survival]
    \label{thm:A-positive-zero-survival}
    Let $0<x_1<\cdots<x_N$ and $d_A<2$. Put
        \begin{equation*}
            \nu_0=1-\frac{d_A}{2},\qquad p_0=N\nu_0,\qquad D_0=\kappa\delta+4\nu_0=4-d_A+\kappa(N-1)>0,
        \end{equation*}
    and let $g(x)=(\prod_i x_i)^{1/N}$ and $\overline{x}=N^{-1}\sum_i x_i$. Then, for every $t>0$,
        \begin{equation}
            \label{eq:A-positive-zero-bounds}
            \left(\frac{g(x)}{\overline{x}+D_0t}\right)^{p_0} \le\mathbb P_x(T_0>t) \le F_{p_0}(t,Ng(x)).
        \end{equation}
    In particular,
        \begin{equation*}
            \mathbb P_x(T_0>t)\asymp t^{-N(1-d_A/2)} \qquad(t\to\infty).
        \end{equation*}
    The process is continued through positive collisions in the definition of $T_0$.
\end{theorem}

For $\kappa=1$ and $0<d_A<2$, the power $N(1-d_A/2)$ already follows from the real Wishart formula in \cite[Proposition~2.6]{MR2105711}. The complex Laguerre formula in \cite[Proposition~7]{Demni2007} gives the same power for $\kappa=2$ and $0<d_A<2$. For $\kappa\ge1$ and $1\le d_A<2$, the first-zero time also falls within the type-$B$ chamber-exit setting of \cite{MR2470522}. Here we obtain explicit two-sided bounds by scalar comparison and a bounded submartingale argument. They hold for every $\kappa>0$ and every $d_A<2$, including $d_A<0$, and allow positive collisions before the first zero. The conclusion is a matching polynomial order; it does not assert an asymptotic equivalent with a specified constant.

Finally, consider the first pair collision from an arbitrary strictly ordered point
    \begin{equation*}
        x\in\mathcal C_A^\circ =\{x\in\mathbb R^N:x_1<\cdots<x_N\}.
    \end{equation*}
The coordinates may have different signs, and one may be zero. Set
    \begin{equation}
        \tau_A=\inf\{t\ge0:X_i(t)=X_{i+1}(t) \text{ for some }1\le i<N\},
    \end{equation}
and write
    \begin{equation*}
        \Delta(x)=\prod_{i<j}(x_j-x_i),\qquad S(x)=\sum_i|x_i|,\qquad M_N=\binom N2.
    \end{equation*}
An individual coordinate hitting or crossing zero does not stop $\tau_A$.

\begin{theorem}[First collision and survival bounds]
    \label{thm:A-collision-main}
    Let $N\ge2$, $\delta\in\mathbb R$ and $x\in\mathcal C_A^\circ$. If $\kappa\ge1$, then $\mathbb P_x(\tau_A=\infty)=1$. If $0<\kappa<1$, put
        \begin{equation*}
            q=1-\kappa,\qquad p=M_Nq,\qquad K=\kappa N|\delta|+2(2-\kappa)M_N.
        \end{equation*}
    Then $\tau_A<\infty$ almost surely and, for every $t>0$,
        \begin{equation}
            \label{eq:A-collision-main-bounds}
            \frac{\Delta(x)^q}{[S(x)+Kt]^p} \le\mathbb P_x(\tau_A>t) \le\prod_{i=1}^{N-1}F_q(t,x_{i+1}-x_i).
        \end{equation}
    Consequently, for some $c_x,C_x>0$ and all sufficiently large $t$,
        \begin{equation}
            \label{eq:A-collision-powers}
            c_x t^{-M_N(1-\kappa)} \le\mathbb P_x(\tau_A>t) \le C_x t^{-(N-1)(1-\kappa)}.
        \end{equation}
    The powers agree for $N=2$.
\end{theorem}

The collision threshold $\kappa=1$ is consistent with the earlier Brownian, Dunkl and nonnegative beta-Wishart results \cite{CepaLepingle1997,Demni2009,JourdainKahn2022}. The statement here covers the whole signed chamber, every real $\delta$, and paths which have already crossed zero. The lower bound uses a power of the Vandermonde function. The upper bound is obtained from all adjacent gaps. For $N>2$ the general bounds need not have matching powers. For $0<\kappa<1$, Lemma~\ref{lem:block}
gives first-collision survival of order $t^{-\binom n2(1-\kappa)}$ for an autonomous nonnegative system of $n\ge2$ particles with constant drift $\kappa d$, $d\ge n-1$, and initial state $0\le v_1<\cdots<v_n$. Its overlap with the regime where zero is not hit is consistent with the type-$B$ chamber-exit formulas in \cite{MR2470522}; our estimate also allows visits to zero before the first collision.

We also study collisions after the first one and from colliding initial states. For every $\kappa>0$, three or more particles cannot collide at a nonzero location at any positive time. For $\kappa\ge1$, Theorem~\ref{thm:two_particle_collisions} excludes pair collisions at nonzero locations as well. Away from zero, this agrees with the collision mechanisms in \cite{CepaLepingle2007,AndrausHufnagelMalecki2025}. At zero the conclusion for three particles can fail. For example, when $N=3$, $\delta=0$ and $x=(-a,0,b)$ with $a,b>0$, the solution has the form $(-U,0,V)$, where $U$ and $V$ are independent squared Bessel processes of dimension $2\kappa$. Their sum has dimension $4\kappa$, so a triple collision at zero occurs in finite time almost surely exactly when $\kappa<1/2$. This example shows why collisions at zero must be treated separately.

The next section constructs the solutions and studies nonuniqueness of the indicator equation. We then study the zero boundary and prove the first-zero survival bounds. The following section treats local collisions, and the final section proves the first-collision bounds.

\section{Construction of solutions}

\subsection{Existence for the indicator equation}

To apply the general existence theorem, we write the system in root coordinates. Let $e_1,\ldots,e_N$ be the standard basis of $\mathbb R^N$, and set
    \begin{equation*}
        R_+=\{e_j-e_i:1\le i<j\le N\}.
    \end{equation*}
For $\alpha=e_j-e_i\in R_+$, define
    \begin{equation*}
        k_\alpha(x)=|x_i|+|x_j|.
    \end{equation*}
Since $\langle x,e_j-e_i\rangle=x_j-x_i$, equation \eqref{eq:A-particle-genuine} has the equivalent root form
    \begin{equation}
        \label{eq:A-root-genuine}
        dX_i(t) = 2\sqrt{|X_i(t)|}\,dB_i(t) + \kappa\left(\delta+ \sum_{\alpha\in R_+}\frac{k_\alpha(X(t))\alpha_i}{\langle X(t),\alpha\rangle} \right)dt.
    \end{equation}
Similarly, the indicator equation \eqref{eq:A-particle-indicator} becomes
    \begin{equation}
        \label{eq:A-root-indicator}
        dX_i(t) = 2\sqrt{|X_i(t)|}\,dB_i(t) + \kappa\left(\delta+\sum_{\alpha\in R_+}\frac{k_\alpha(X(t))\alpha_i}{\langle X(t),\alpha\rangle} \mathbf 1_{\{\langle X(t),\alpha\rangle>0\}}\right)dt.
\end{equation}
We apply \cite[Theorem~3]{Malecki2026}. In the notation of that theorem,
    \begin{equation*}
        \sigma(x)=\operatorname{diag}(2\sqrt{|x_1|},\ldots,2\sqrt{|x_N|}), \qquad b(x)=\kappa\delta(1,\ldots,1),
    \end{equation*}
and the root coefficient attached to $e_j-e_i$ is $\kappa k_{e_j-e_i}(x)= \kappa(|x_i|+|x_j|)$. The coefficients are continuous, and
    \begin{equation*}
        |\sigma_i^2(x)-\sigma_j^2(x)| = 4\bigl||x_i|-|x_j|\bigr| \le4|x_i-x_j|.
    \end{equation*}
Thus the diffusion compatibility condition holds. Moreover,
    \begin{equation*}
        \kappa k_{e_j-e_i}(x) = \frac{\kappa}{4} \bigl(\sigma_i^2(x)+\sigma_j^2(x)\bigr),
    \end{equation*}
which verifies the normal dominance condition $(D1)$. On every collision block, the regular drifts are equal. The interactions with particles outside the block are also equal for all coordinates in that block. Every normal detector therefore has zero regular drift on the face, verifying $(D2)$. Hence \cite[Theorem~3]{Malecki2026} gives a weak solution of \eqref{eq:A-particle-indicator} in $\mathcal C_A$, with finite singular drift integrals on compact time intervals.

To obtain global existence, we verify the radial growth condition $(G1)$ of
\cite{Malecki2026}. We have
    \begin{equation*}
        \langle x,b(x)\rangle\le C(1+|x|^2),\qquad \frac12\operatorname{Tr}(\sigma\sigma^\top)(x)\le C(1+|x|^2),
    \end{equation*}
because both terms are bounded by a constant times $1+\sum_i |x_i|$. Moreover,
    \begin{equation*}
        \sum_{\alpha\in R_+}\kappa k_\alpha(x) = \kappa\sum_{1\le i<j\le N}(|x_i|+|x_j|) = \kappa(N-1)\sum_{i=1}^N|x_i| \le C(1+|x|^2).
    \end{equation*}
Thus $(G1)$ holds, and \cite[Theorem~5]{Malecki2026} implies non-explosion. Consequently we obtain the following.

\begin{proposition}
    \label{prop:A-indicator-existence}
    Let $N\ge2$, $\kappa>0$, and $\delta\in\mathbb R$. For every $x\in\mathcal C_A$, there exists a global weak solution $X$ of
        \begin{equation*}
            dX_i(t) = 2\sqrt{|X_i(t)|}\,dB_i(t) + \kappa\delta\,dt + \kappa\sum_{j\ne i} \frac{|X_i(t)|+|X_j(t)|}{X_i(t)-X_j(t)} \mathbf 1_{\{X_i(t)\ne X_j(t)\}}\,dt,
        \end{equation*}
    with $X(t)\in\mathcal C_A$ for all $t\ge0$. Moreover, the singular drift integrals are finite on compact time intervals.
\end{proposition}

\subsection{Zero occupation and the non-resonant case}
We have already obtained weak existence for the indicator equation. We still need to know whether its solutions are non-sticky. In this subsection we identify the configurations where collisions may have positive occupation time, and then prove pathwise uniqueness in the non-sticky class. These two steps will give strong existence outside resonance.

For $x\in\mathcal C_A$, write
    \begin{equation*}
        \ell(x)=\#\{i:x_i<0\},\qquad m(x)=\#\{i:x_i=0\},\qquad n(x)=\#\{i:x_i>0\}.
    \end{equation*}
Define
    \begin{equation*}
        \mathcal R_N = \{n-\ell: \ell,m,n\in\mathbb Z_{\ge0},\ \ell+m+n=N,\ m\ge2\} = \mathbb Z\cap[2-N,N-2]
    \end{equation*}
and
    \begin{equation*}
        \mathcal Z_\delta = \{x\in\mathcal C_A: m(x)\ge2,\ \delta=n(x)-\ell(x)\}.
    \end{equation*}
The set $\mathcal R_N$ contains the values of $\delta$ for which the common indicator drift can vanish at a zero block of size at least two. For a fixed $\delta$, the corresponding configurations form $\mathcal Z_\delta$. The following calculation explains these definitions.

\begin{lemma}[Effective drift at zero]
    \label{lem:zero-block-drift}
    At a configuration $x\in\mathcal C_A$ with $x_i=0$, the drift of the $i$-th coordinate in the indicator equation is
        \begin{equation*}
            \kappa\bigl(\delta+\ell(x)-n(x)\bigr).
        \end{equation*}
    For a zero block of size $m$, with $\ell$ negative and $n=N-\ell-m$ positive coordinates, this equals $\kappa(\delta+2\ell+m-N)$.
\end{lemma}

\begin{proof}
    An interaction with a negative coordinate contributes $+1$, an interaction with a positive coordinate contributes $-1$, and interactions with other zero coordinates vanish by the indicators.
\end{proof}

The next lemma shows that pairwise collisions can have positive occupation time only on $\mathcal Z_\delta$. At zero, we use Tanaka's formula and the vanishing of both one-sided local times to obtain a condition on the indicator drift. Nonzero collisions are treated separately by applying the occupation-time formula to the gaps.

\begin{lemma}[Zero local time and possible sticking]
    \label{lem:A-zero-occupation}
    Let $X$ be a global solution of \eqref{eq:A-particle-indicator}. Both one-sided local times of every $X_i$ at zero vanish. Moreover, almost surely,
        \begin{equation}
            \label{eq:A-zero-drift-support}
            \mathbf 1_{\{X_i(t)=0\}} \bigl(\delta+\ell(X(t))-n(X(t))\bigr)=0 \quad\text{for Lebesgue-a.e. }t.
        \end{equation}
    For every $T>0$ and every $i<j$,
        \begin{equation}
            \label{eq:A-collisions-outside-resonance}
            \int_0^T \mathbf 1_{\{X_i(t)=X_j(t)\}} \mathbf 1_{\{X(t)\notin\mathcal Z_\delta\}}\,dt=0 \quad\text{a.s.}
        \end{equation}
\end{lemma}

\begin{proof}
    We have
        \begin{equation*}
            d\langle X_i\rangle_t=4|X_i(t)|\,dt.
        \end{equation*}
    Write $L_T^a(X_i)$ for right local time. The occupation-time formula gives
        \begin{equation*}
            \int_0^1\frac{L_T^a(X_i)}{a}\,da = 4\int_0^T\mathbf 1_{\{0<X_i(t)<1\}}\,dt <\infty.
        \end{equation*}
    Right continuity in the level implies $L_T^0(X_i)=0$. The same argument on negative levels, using the left limit in the level, gives $L_T^{0-}(X_i)=0$. Taking rational $T$ and using continuity in time gives these conclusions simultaneously for all $T$. The two Tanaka formulas consequently imply
        \begin{equation*}
            \int_0^t\mathbf 1_{\{X_i(s)=0\}}\,dX_i(s)=0, \qquad t\ge0.
        \end{equation*}
    The martingale part of this integral vanishes. Lemma~\ref{lem:zero-block-drift} therefore gives
        \begin{equation*}
            \int_0^t \mathbf 1_{\{X_i(s)=0\}} \kappa\bigl(\delta+\ell(X(s))-n(X(s))\bigr)\,ds=0.
        \end{equation*}
    Differentiating this absolutely continuous identity proves \eqref{eq:A-zero-drift-support}.

    For $D_{ij}=X_j-X_i$, independence of the Brownian coordinates gives
        \begin{equation*}
            d\langle D_{ij}\rangle_t = 4\bigl(|X_i(t)|+|X_j(t)|\bigr)\,dt.
        \end{equation*}
    The occupation-time formula at the single level zero yields
        \begin{equation*}
            \int_0^T \mathbf 1_{\{X_i(t)=X_j(t)\}} \bigl(|X_i(t)|+|X_j(t)|\bigr)\,dt=0.
        \end{equation*}
    Thus nonzero collisions have zero occupation time. At a zero collision, $m(X(t))\ge2$, and \eqref{eq:A-zero-drift-support} excludes occupation unless $\delta=n(X(t))-\ell(X(t))$. This proves \eqref{eq:A-collisions-outside-resonance}.
\end{proof}

The uniqueness argument below applies to every $\delta$, also at resonance. We compare two non-sticky solutions with the same initial condition and driven by the same Brownian motion. Tanaka's formula reduces the question to the sign of the interaction terms, which is controlled by the monotonicity of the pair interaction.

\begin{lemma}[Pathwise uniqueness in the non-sticky class]
    \label{lem:A-nonsticky-uniqueness}
    Two non-sticky solutions of \eqref{eq:A-particle-indicator}, with the same initial condition and the same Brownian motion, are indistinguishable.
\end{lemma}

\begin{proof}
    Let $X,\widehat X$ be such solutions and put $D_i=X_i-\widehat X_i$. Since
        \begin{equation*}
            |2\sqrt{|u|}-2\sqrt{|v|}|^2\le4|u-v|,
        \end{equation*}
    we have $d\langle D_i\rangle_t\le4|D_i(t)|\,dt$. The local-time argument in Lemma~\ref{lem:A-zero-occupation} shows that both one-sided local times of $D_i$ at zero vanish. For $u>v$, set
        \begin{equation*}
            \Psi(u,v)=\frac{|u|+|v|}{u-v}.
        \end{equation*}
    This function is non-increasing in $u$ and non-decreasing in $v$. At times when both configurations have distinct coordinates, write $s_i=\operatorname{sgn}(D_i)$. The contribution of a pair $i<j$ to the drift in Tanaka's formula for $\sum_i|D_i|$ is
        \begin{equation*}
            \kappa(s_j-s_i)\bigl(\Psi(X_j,X_i)-\Psi(\widehat X_j,\widehat X_i)\bigr) \le0.
        \end{equation*}
    The expression vanishes if $s_i=s_j$; otherwise its sign follows from the two monotonicities. Non-stickiness makes the exceptional collision times a Lebesgue-null set. The constant drifts cancel. Consequently, $\sum_i|D_i|$ is a nonnegative local supermartingale started from zero. After localization and taking expectations, it follows that it vanishes identically.
\end{proof}

We use the Yamada-Watanabe principle within the class of non-sticky solutions with locally integrable singular drifts. These are measurable path properties and are preserved by the conditional coupling in its proof. Thus weak existence and pathwise uniqueness within this class yield strong existence and uniqueness in law within the same class. Outside resonance, $\mathcal Z_\delta$ is empty, so every indicator solution is non-sticky. We can therefore combine the weak existence result with the preceding uniqueness argument to obtain the following theorem.
\begin{theorem}[Non-resonant solutions]
    \label{thm:A-nonresonant-genuine}
    Let $N\ge2$, $\kappa>0$, and $\delta\notin\mathcal R_N$. For every $x\in\mathcal C_A$, equation \eqref{eq:A-particle-genuine} has a global strong solution. Pathwise uniqueness holds among all solutions of the genuine equation. Every solution of the indicator equation is non-sticky.
\end{theorem}

\begin{proof}
    Proposition~\ref{prop:A-indicator-existence} gives a global weak indicator solution. Since $\mathcal Z_\delta=\varnothing$, Lemma~\ref{lem:A-zero-occupation} makes every indicator solution non-sticky. The indicators can therefore be removed. Lemma~\ref{lem:A-nonsticky-uniqueness} gives pathwise uniqueness, and the Yamada-Watanabe principle gives strong existence.
\end{proof}

\begin{rem}[The role of resonance]
    At a resonant zero block the common indicator drift is zero. Lemma~\ref{lem:A-zero-occupation} therefore does not exclude occupation of that block. Below we construct both a non-sticky continuation and a solution whose designated zero block stays pinned.
\end{rem}

\subsection{One-sided squared Bessel systems}

To construct the continuations below, we shall use smaller systems which remain nonnegative. Their coordinates may start at zero, but must spend zero Lebesgue time there. We first prove nonnegativity for the indicator equation and then obtain the required strong solutions. For $r\ge0$, write
    \begin{equation*}
        \mathcal C_r^+=\{0\le y_1\le\cdots\le y_r\}.
    \end{equation*}
The first lemma applies to every indicator solution, without assuming non-stickiness. We use the sum of the negative parts of all coordinates, because the interactions between negative coordinates cancel in its equation.

\begin{lemma}[Invariance of the nonnegative chamber]
    \label{lem:A-indicator-nonnegative}
    Suppose that $\delta\ge N-1$. Every solution of \eqref{eq:A-particle-indicator} started in $\{0\le x_1\le\cdots\le x_N\}$ remains in this set.
\end{lemma}

\begin{proof}
    Put
        \begin{equation*}
            S(t)=\sum_{i=1}^N X_i(t)^-, \qquad K(t)=\#\{i:X_i(t)<0\}.
        \end{equation*}
    By Lemma~\ref{lem:A-zero-occupation}, the local-time terms in Tanaka's formula vanish. Interactions within the negative group cancel in pairs. Each negative coordinate has $N-K(t)$ interactions with nonnegative coordinates, each contributing $-1$. Hence
        \begin{equation*}
            dS(t) = -2\sum_{i=1}^N \mathbf 1_{\{X_i(t)<0\}}\sqrt{|X_i(t)|}\,dB_i(t) \quad -\kappa K(t)\bigl(\delta-N+K(t)\bigr)\,dt.
        \end{equation*}
    For $K(t)\ge1$,
        \begin{equation*}
            \delta-N+K(t)\ge K(t)-1\ge0.
        \end{equation*}
    For $K(t)=0$, the drift vanishes. Thus $S$ is a nonnegative local supermartingale started from zero, and consequently $S(t)=0$ for every $t\ge0$.
\end{proof}
For an $r$-particle system, the preceding lemma gives nonnegativity under $d\ge r-1$. We shall work with $d\ge r$, so that every zero coordinate has strictly positive indicator drift, also when several coordinates are zero together. Equation~\eqref{eq:A-zero-drift-support} then rules out positive time spent at the origin.

For an integer $r\ge0$ and $d\ge r$, consider
    \begin{equation}
        \label{eq:one-sided-system}
        dY_a=2\sqrt{Y_a}\,dW_a+ \kappa\left(d+\sum_{b\ne a} \frac{Y_a+Y_b}{Y_a-Y_b}\right)dt, \qquad a=1,\ldots,r.
    \end{equation}
Here $W_1,\ldots,W_r$ are independent standard Brownian motions. For $r=0$, this is the empty system.

\begin{lemma}[One-sided auxiliary systems]
    \label{lem:one-sided-BESQ-system}
    For every $y\in\mathcal C_r^+$, equation \eqref{eq:one-sided-system} has a pathwise unique global strong solution in $\mathcal C_r^+$, with locally absolutely integrable singular drifts. The solution is pairwise non-sticky, and, almost surely,
        \begin{equation*}
            \int_0^T\mathbf1_{\{Y_a(t)=0\}}\,dt=0, \qquad a=1,\ldots,r,\quad T>0.
        \end{equation*}
\end{lemma}

\begin{proof}
    The case $r=0$ is empty. For $r=1$, the equation is the classical squared Bessel equation of dimension $\kappa d>0$. Its global strong solution is nonnegative and spends zero Lebesgue time at zero. Let $r\ge2$. Apply Proposition~\ref{prop:A-indicator-existence} with $N=r$ and $\delta=d$, initially using the absolute-value coefficients. Since $d\ge r>r-1$, Lemma~\ref{lem:A-indicator-nonnegative} shows that the resulting global indicator solution stays in $\mathcal C_r^+$.

    At a time when $Y_a=0$, let $s\ge1$ be the number of zero coordinates. The drift of $Y_a$ at that time is
        \begin{equation*}
            \kappa(d-r+s)\ge\kappa s>0.
        \end{equation*}
    Equation~\eqref{eq:A-zero-drift-support}, applied to this $r$-particle system, implies that every coordinate spends zero Lebesgue time at zero. Nonzero collisions also have zero occupation time by Lemma~\ref{lem:A-zero-occupation}. Thus the solution is pairwise non-sticky, and removing the indicators gives \eqref{eq:one-sided-system}.

    The singular drift integrability and non-explosion follow from the indicator existence proposition. Lemma~\ref{lem:A-nonsticky-uniqueness} gives pathwise uniqueness, and the Yamada-Watanabe principle gives strong existence.
\end{proof}
The auxiliary systems may still visit zero or have pair collisions. What we use below is that both individual zero times and pairwise collision times have zero Lebesgue measure. The assertion about individual coordinates is needed in addition to pairwise non-stickiness, which by itself does not exclude a single coordinate remaining at zero.

\subsection{Balanced continuation from a resonant zero face}

We construct a non-sticky continuation from a designated resonant zero block, allowing additional zeros in the side groups. The construction splits the designated block equally between the two half-lines; if its size is odd, the middle particle remains at zero. The equal split preserves the difference $n-\ell$ between the sizes of the side groups. A single coordinate may remain fixed at zero, since non-stickiness concerns pairs of coordinates. We also allow additional zeros in the side groups, as this will be needed at the stopping time used in the global construction.
\begin{proposition}[Balanced resonant continuation]
    \label{prop:balanced-resonant-continuation}
    Let $N\ge2$, $\kappa>0$, and let
        \begin{equation*}
            \ell+m+n=N,\qquad m\ge2,\qquad \delta=n-\ell.
        \end{equation*}
    Suppose that $x\in\mathcal C_A$ satisfies
        \begin{equation*}
            x_1\le\cdots\le x_\ell\le0 = x_{\ell+1}=\cdots=x_{\ell+m} \le x_{\ell+m+1}\le\cdots\le x_N,
        \end{equation*}
    with the obvious conventions when $\ell=0$ or $n=0$. The designated zero block need not be maximal. Then there exists a global non-sticky weak solution of \eqref{eq:A-particle-genuine} started from $x$. More precisely, put
        \begin{equation*}
            h=\left\lfloor\frac m2\right\rfloor,\qquad q=m-2h\in\{0,1\},\qquad L=\ell+h,\qquad R=n+h.
        \end{equation*}
    The solution can be chosen so that
        \begin{equation*}
            X_1(t)\le\cdots\le X_L(t)\le0 \le X_{L+q+1}(t)\le\cdots\le X_N(t)
        \end{equation*}
    for all $t\ge0$, and, when $q=1$,
        \begin{equation*}
            X_{L+1}(t)=0,\qquad t\ge0.
        \end{equation*}
\end{proposition}

\begin{proof}
    The first $L$ initial coordinates are non-positive and the last $R$ coordinates are non-negative.  Define
        \begin{equation*}
            u_a=-x_{L+1-a},\qquad a=1,\ldots,L,
        \end{equation*}
    and
        \begin{equation*}
            v_b=x_{L+q+b},\qquad b=1,\ldots,R.
        \end{equation*}
    For particles on opposite sides of zero, the interaction quotient is $+1$ or $-1$ whenever it is defined. These terms can be included in the constant drifts, so the two one-sided systems can be constructed independently. The parameters chosen below account for these interactions and for the possible middle coordinate fixed at zero. Then $u\in\mathcal C_L^+$ and $v\in\mathcal C_R^+$.  By Lemma~\ref{lem:one-sided-BESQ-system}, we may construct independent one-sided systems $U$ and $V$, started respectively from $u$ and $v$, with dimension parameters
        \begin{equation*}
            d_U=L+q,\qquad d_V=R+q.
        \end{equation*}
    Thus
        \begin{align}
            dU_a &= 2\sqrt{U_a}\,dW_a^- + \kappa\left(L+q+ \sum_{c\ne a}\frac{U_a+U_c}{U_a-U_c} \right)dt,
            \label{eq:balanced-U}\\
            dV_b &= 2\sqrt{V_b}\,dW_b^+ + \kappa\left(R+q+\sum_{c\ne b}\frac{V_b+V_c}{V_b-V_c}\right)dt.
            \label{eq:balanced-V}
        \end{align}
    The assumptions of the lemma hold because $L+q\ge L$ and $R+q\ge R$. Define
        \begin{equation*}
            X_i(t)=
            \begin{cases}
                -\,U_{L+1-i}(t),&1\le i\le L,\\[1mm]
                0,&q=1\ \text{and }i=L+1,\\[1mm]
                V_{i-L-q}(t),&L+q+1\le i\le N.
            \end{cases}
        \end{equation*}
    Using $-W^-_{L+1-i}$ for the Brownian motion of the $i$-th negative coordinate, $W^+_{i-L-q}$ for the positive coordinates, and, when $q=1$, an additional independent Brownian motion for the pinned coordinate, we obtain an $N$-dimensional Brownian motion $B$.  Since $U$ and $V$ stay in their one-sided Weyl chambers, $X(t)\in\mathcal C_A$. We verify the equation.  All identities below are understood $dt$-a.e.; the exceptional sets are contained in zero sets of coordinates of $U$ or $V$, or in their internal collision sets. All these sets have zero Lebesgue measure by Lemma~\ref{lem:one-sided-BESQ-system}. If $X_i=V_b>0$, then each of the $L$ non-positive particles and, when $q=1$, the pinned zero particle contributes $+1$. Since
        \begin{equation*}
            \delta+L+q = n-\ell+\ell+h+q = n+h+q = R+q,
        \end{equation*}
    the constant part of the drift agrees with \eqref{eq:balanced-V}; the interactions with the other positive particles agree identically. If $X_i=-U_a<0$, then every positive particle and the possible pinned zero particle contributes $-1$, while for another negative particle $X_j=-U_c$,
        \begin{equation*}
            \frac{|X_i|+|X_j|}{X_i-X_j} = - \frac{U_a+U_c}{U_a-U_c}.
        \end{equation*}
    Moreover,
        \begin{equation*}
            R+q-\delta = n+h+q-(n-\ell) = \ell+h+q = L+q.
        \end{equation*}
    After multiplying \eqref{eq:balanced-U} by $-1$, this is exactly the drift required for $X_i$. It remains only to consider the case $q=1$.  For the pinned coordinate $X_{L+1}=0$, the diffusion coefficient vanishes.  The $L$ particles to its left contribute $+1$ and the $R$ particles to its right contribute $-1$, so its drift is
        \begin{equation*}
            \kappa(\delta+L-R) = \kappa\bigl(n-\ell+\ell+h-(n+h)\bigr) = 0.
        \end{equation*}
    Thus the pinned coordinate also satisfies the equation.

The singular drift integrals internal to the negative and positive groups are finite by Lemma~\ref{lem:one-sided-BESQ-system}.  Every cross-interaction is equal to $+1$ or $-1$ whenever it is defined, and may be assigned any value on its zero-Lebesgue-time collision set.  Hence all cross-interaction integrals are finite.

Finally, collisions inside either side have zero occupation time by the lemma.  A collision between a non-positive and a non-negative particle can occur only when both are zero, and its time set is contained in the zero set of either coordinate.  In the case $q=1$, a collision with the pinned coordinate occurs exactly when the corresponding side coordinate is zero. All such sets have zero Lebesgue measure.  Therefore $X$ is pairwise non-sticky and solves the genuine equation globally.
\end{proof}

\subsection{Sticky solutions of the indicator equation}

Although Proposition~\ref{prop:balanced-resonant-continuation} produces a non-sticky solution, the resonant face can also support a sticky solution of the indicator equation. This explains why the indicators cannot simply be removed from every indicator solution. The indicators allow all interactions inside the pinned zero block to be set to zero. The resonance relation then makes the total drift of each pinned coordinate vanish for almost every time. This allows us to keep the whole designated block at zero.

\begin{proposition}[Sticky resonant solutions]
    \label{prop:sticky-resonant-solutions}
    Let $N\ge2$, $\kappa>0$, and suppose that
        \begin{equation*}
            \delta=n-\ell, \qquad \ell+m+n=N, \qquad m\ge2 ,
        \end{equation*}
    for some integers $\ell,n,m\ge0$. Let $x\in\mathcal C_A$ satisfy
        \begin{equation*}
            x_1\le\cdots\le x_\ell\le0 = x_{\ell+1}=\cdots=x_{\ell+m} \le x_{\ell+m+1}\le\cdots\le x_N ,
        \end{equation*}
    with the obvious convention if $\ell=0$ or $n=0$. Then the indicator equation admits a global strong solution $X$ such that
        \begin{equation*}
            X_{\ell+1}(t)=\cdots=X_{\ell+m}(t)=0, \qquad t\ge0.
        \end{equation*}
    In particular, the resonant zero face is invariant for at least one strong solution of the indicator equation.
\end{proposition}

\begin{proof}
    Set
        \begin{equation*}
            u_a=-x_{\ell+1-a},\qquad a=1,\ldots,\ell, \qquad v_b=x_{\ell+m+b},\qquad b=1,\ldots,n .
        \end{equation*}
    Each side system must now account for interactions with all $m$ pinned particles. This explains the parameters $\ell+m$ and $n+m$ used below. Let $U$ be the one-sided system of Lemma~\ref{lem:one-sided-BESQ-system} with $r=\ell$ and $d=\ell+m$, and let $V$ be an independent one-sided system with $r=n$ and $d=n+m$, started from $u$ and $v$, respectively. Define
        \begin{equation*}
            X_i(t)=
            \begin{cases}
                -\,U_{\ell+1-i}(t), & 1\le i\le \ell,\\[1mm]
                0, & \ell+1\le i\le \ell+m,\\[1mm]
                V_{i-\ell-m}(t), & \ell+m+1\le i\le N.
            \end{cases}
        \end{equation*}
    Together with independent Brownian motions for the pinned zero coordinates, this defines an $N$-dimensional Brownian motion $B$. Since $U,V$ remain in their one-sided Weyl chambers, $X(t)\in\mathcal C_A$.

    We check the drift. By Lemma~\ref{lem:one-sided-BESQ-system}, the side particles spend zero Lebesgue time at $0$, so the following identities hold $dt$-a.e. For $X_i=V_b>0$, each particle on the negative side and in the pinned zero block contributes $+1$. Hence
        \begin{equation*}
            \delta+\ell+m = n-\ell+\ell+m = n+m,
        \end{equation*}
    which is exactly the dimension parameter of the $V$-system. For $X_i=-U_a<0$, each pinned zero particle and each positive particle contributes $-1$, while the interactions with the other negative particles are the negatives of the corresponding $U$-interactions:
        \begin{equation*}
            \frac{|X_i|+|X_j|}{X_i-X_j} = - \frac{U_a+U_c}{U_a-U_c}.
        \end{equation*}
    Thus $\delta-(m+n)=n-\ell-m-n=-(\ell+m)$,  so the drift agrees with that of $-U_a$.

    Finally, for a pinned particle $X_i=0$, the diffusion coefficient is zero, the interactions inside the pinned block vanish by the indicators, and the remaining cross-interactions give $\ell$ contributions $+1$ and $n$ contributions $-1$. Therefore its drift is $\kappa(\delta+\ell-n) = 0$. Thus every coordinate satisfies the indicator equation.
    The finiteness of the singular drift integrals follows from
    the corresponding property of $U,V$, since all cross-interaction
    terms are bounded by $1$ in absolute value, apart from sets
    killed by the indicators.
    Moreover, since $U$ and $V$ are strong solutions driven by the
    corresponding coordinates of $B$, the constructed process $X$
    is adapted to the filtration generated by $B$.
    Hence $X$ is a global strong solution of the indicator equation
    with the designated zero block pinned at zero.
\end{proof}

\subsection{Global existence for all parameters}

    We now use the balanced construction to obtain a non-sticky solution from any initial point. We keep a weak indicator solution up to a suitable stopping time and use the balanced continuation afterwards. We choose a closed stopping set containing $\mathcal Z_\delta$, so that the state at a finite first hitting time belongs to this set. Define the finite set of block patterns
        \begin{equation*}
            \mathfrak T_\delta = \{(\ell,m,n)\in\mathbb Z_{\ge0}^3: \ell+m+n=N,\ m\ge2,\ \delta=n-\ell\}.
        \end{equation*}
    For $(\ell,m,n)\in\mathfrak T_\delta$, set
        \begin{equation*}
            F_{\ell,m,n} = \{x\in\mathcal C_A: x_{\ell+1}=\cdots=x_{\ell+m}=0\}, \qquad \widehat{\mathcal Z}_\delta = \bigcup_{(\ell,m,n)\in\mathfrak T_\delta}F_{\ell,m,n}.
        \end{equation*}
    The union is empty if $\mathfrak T_\delta=\varnothing$. It is a finite union of closed sets and therefore closed. Moreover, $\mathcal Z_\delta\subseteq\widehat{\mathcal Z}_\delta$.

    The indices in a selected pattern need not equal the actual sign counts of a point in $F_{\ell,m,n}$, since additional side coordinates may also be zero. The larger set $\widehat{\mathcal Z}_\delta$ may be reached before $\mathcal Z_\delta$. This causes no difficulty, because the balanced construction allows additional zeros in the side groups and can already be used at that time. The next lemma shows that the original solution is non-sticky up to this stopping time.

\begin{lemma}[Non-stickiness before the gluing time]
    \label{lem:nonsticking-before-resonance}
    Let $X$ be a global solution of \eqref{eq:A-particle-indicator}, and set
        \begin{equation*}
            \tau = \inf\{t\ge0:X(t)\in\widehat{\mathcal Z}_\delta\}.
        \end{equation*}
    Then, almost surely, for every $T>0$ and $i<j$,
        \begin{equation*}
            \int_0^{T\wedge\tau} \mathbf 1_{\{X_i(t)=X_j(t)\}}\,dt=0.
        \end{equation*}
    On $\{\tau<\infty\}$, $X(\tau)\in\widehat{\mathcal Z}_\delta$.
\end{lemma}

\begin{proof}
    Before $\tau$, the process lies outside $\mathcal Z_\delta$. Equation~\eqref{eq:A-collisions-outside-resonance} gives the occupation assertion; the endpoint contributes nothing. The stopping-time property and the assertion about $X(\tau)$ follow from continuity and closedness.
\end{proof}

    We need at most one switch, since the balanced continuation is already global and non-sticky. The construction below first gives weak existence in the non-sticky class. Strong existence will then follow from pathwise uniqueness and the Yamada-Watanabe principle.

\begin{proof}[Proof of Theorem~\ref{thm:A-genuine-existence-uniqueness}]
    By Proposition~\ref{prop:A-indicator-existence}, let $X^0$ be a global indicator solution driven by an $N$-dimensional Brownian motion $B$. Set
    \begin{equation*}
        \tau = \inf\{t\ge0:X^0(t)\in\widehat{\mathcal Z}_\delta\}.
    \end{equation*}
    Retain $X^0$ up to $\tau$. Lemma~\ref{lem:nonsticking-before-resonance} shows that it is non-sticky on this interval. If $\tau=\infty$, no continuation is needed.

    On $\{\tau<\infty\}$, choose the first triple $(\ell,m,n)$, in a fixed ordering of $\mathfrak T_\delta$, such that $X^0(\tau)\in F_{\ell,m,n}$. This choice is $\mathcal F_\tau$-measurable. Put
        \begin{equation*}
            h=\lfloor m/2\rfloor,\qquad q=m-2h,\qquad L=\ell+h,\qquad R=n+h.
        \end{equation*}
    The state $X^0(\tau)$ satisfies the hypotheses of Proposition~\ref{prop:balanced-resonant-continuation}, including the permitted additional zeros in the side groups. Use the shifted Brownian coordinates
        \begin{equation*}
            W_a^-(t) = -\bigl(B_{L+1-a}(\tau+t)-B_{L+1-a}(\tau)\bigr), \qquad a=1,\ldots,L,
        \end{equation*}
    and
        \begin{equation*}
            W_b^+(t) = B_{L+q+b}(\tau+t)-B_{L+q+b}(\tau), \qquad b=1,\ldots,R.
        \end{equation*}
    For each selected pattern, these are independent Brownian motions after $\tau$, independent of the past. By Lemma~\ref{lem:one-sided-BESQ-system}, solve the one-sided equations with parameters $L+q$ and $R+q$, from
        \begin{equation*}
            U_a(0)=-X^0_{L+1-a}(\tau), \qquad V_b(0)=X^0_{L+q+b}(\tau).
        \end{equation*}
    Strong existence and pathwise uniqueness allow these solutions to be chosen measurably for the random initial states.

    Reflect $U$, retain the middle coordinate at zero when $q=1$, and append $V$, as in the balanced-continuation proposition. Denote the resulting continuation by $Y$. Define
        \begin{equation*}
            X(t) =
            \begin{cases}
                X^0(t),&t\le\tau,\\
                Y(t-\tau),&t>\tau.
            \end{cases}
        \end{equation*}
    The process is continuous and adapted. Before $\tau$, the indicator equation agrees with the genuine equation because the collision sets have zero occupation time. After $\tau$, the balanced-continuation proposition gives the genuine equation and non-stickiness. Splitting the stochastic and finite-variation integrals at $\tau$ proves the integral equation globally, with the same Brownian motion $B$. The singular integrals are finite on compact intervals on both sides of the gluing time, hence also across $\tau$. All component processes are global. This proves global weak existence in the non-sticky class.

    By the solution convention in the Introduction, every solution of \eqref{eq:A-particle-genuine} is a non-sticky indicator solution. Lemma~\ref{lem:A-nonsticky-uniqueness} therefore gives pathwise uniqueness among all solutions of the genuine equation. The Yamada-Watanabe principle in this class yields global strong existence and uniqueness in law.
\end{proof}

\subsection{Delayed release and nonuniqueness}

At a resonant zero block, we have constructed both a non-sticky continuation and a solution with the block pinned for all time. We can also keep the block at zero up to a prescribed deterministic time and then use the non-sticky solution from the state reached at that time. The next proposition shows that different holding times give different solution laws.

\begin{proposition}[Arbitrary holding times at a resonant zero block]
    \label{prop:A-delayed-release}
    Assume the hypotheses of Proposition~\ref{prop:sticky-resonant-solutions}, and put $I=\{\ell+1,\ldots,\ell+m\}$. Let $B$ be an $N$-dimensional Brownian motion. For every deterministic $a\in[0,\infty]$, the indicator equation \eqref{eq:A-particle-indicator} has a global strong solution $X^{[a]}$, started from $x$ and driven by $B$, such that
        \begin{equation*}
            X_i^{[a]}(t)=0, \qquad i\in I,\quad 0\le t\le a.
        \end{equation*}
    If $a<\infty$, the solution is pairwise non-sticky on $[a,\infty)$. In particular, almost surely, for every $T>0$ and distinct $i,j\in I$,
        \begin{equation*}
            \int_0^T\mathbf 1_{\{X_i^{[a]}(t)=X_j^{[a]}(t)\}}\,dt = T\wedge a,
        \end{equation*}
    where $T\wedge\infty=T$. The solutions corresponding to distinct values of $a$ have distinct laws. Consequently, the indicator equation has uncountably many distinct solution laws from this initial state.
\end{proposition}

\begin{proof}
    The construction in Proposition~\ref{prop:sticky-resonant-solutions} uses strong one-sided solutions. We can therefore carry it out with the prescribed Brownian motion $B$, using its signed and reordered coordinates for the side systems. The Brownian coordinates corresponding to the pinned particles are unused at this stage. Denote the resulting solution by $X^{\rm pin}$.

    For $a=\infty$, set $X^{[\infty]}=X^{\rm pin}$. Let now $a<\infty$. Run $X^{\rm pin}$ until time $a$. By Theorem~\ref{thm:A-genuine-existence-uniqueness} and its extension to random initial conditions, let $Y$ be the non-sticky strong solution started from $X^{\rm pin}(a)$ and driven by $B(a+\cdot)-B(a)$. This use of a random initial state is justified by strong existence and pathwise uniqueness; the state is $\mathcal F_a$-measurable and the future Brownian increments are independent of $\mathcal F_a$.

    Define
        \begin{equation*}
            X^{[a]}(t)=
                \begin{cases}
                    X^{\rm pin}(t),&t\le a,\\
                    Y(t-a),&t>a.
                \end{cases}
        \end{equation*}
    Both pieces solve the indicator equation and their states agree at time $a$. The process is continuous and adapted, and its singular drift integrals are finite on compact intervals. Splitting the stochastic and drift integrals at $a$ shows that $X^{[a]}$ solves the indicator equation with the same Brownian motion $B$.

    For distinct $i,j\in I$, the two coordinates agree throughout $[0,a]$. After time $a$, their collision set has zero Lebesgue measure. This proves the occupation identity. To distinguish the laws, fix $a<b$ and choose $T>a$. For $X^{[a]}$, the occupation integral equals $a$, whereas for $X^{[b]}$ it equals $T\wedge b>a$. This integral is a measurable function of the solution path, so the two laws are different. Since there are uncountably many choices of $a$, the last assertion follows.
\end{proof}

The choices $a=0$ and $a=\infty$ give the non-sticky and permanently pinned solutions, respectively. For a finite positive $a$, the designated block is held at zero until $a$, and the solution is pairwise non-sticky afterwards.

We can now summarize when the indicator equation is pathwise unique for every initial configuration.

\begin{cor}[Uniqueness for all initial configurations]
    \label{cor:A-indicator-nonuniqueness}
    Let $N\ge2$, $\kappa>0$ and $\delta\in\mathbb R$. The indicator equation \eqref{eq:A-particle-indicator} has a pathwise unique global strong solution for every $x\in\mathcal C_A$ if and only if $\delta\notin\mathcal R_N$.
\end{cor}

\begin{proof}
    If $\delta\notin\mathcal R_N$, the assertion follows from Theorem~\ref{thm:A-nonresonant-genuine}. If $\delta\in\mathcal R_N$, choose a corresponding block pattern and an initial state satisfying the assumptions of Proposition~\ref{prop:A-delayed-release}. The choices $a=0$ and $a=\infty$ give distinct strong solutions from this state, driven by the same Brownian motion. Thus pathwise uniqueness does not hold for every initial configuration.
\end{proof}

\section{The zero boundary from the positive chamber}

We now study whether a solution started in the strictly positive chamber can hit zero and, if it does, whether it subsequently enters the negative half-line.  Throughout this section, $X$ denotes the global non-sticky strong solution of \eqref{eq:A-particle-genuine} started from $0<x_1<\cdots<x_N$. We first study the first zero through the product of the coordinates, whose equation can be reduced to a squared Bessel equation by a time change. The question of entering the negative half-line will then be treated separately, using the leftmost particle of the zero block.

\subsection{A random time-change lemma}

The following elementary lemma will be used in the sequel to reduce a one-dimensional semimartingale with random diffusion speed to squared-Bessel form.

Let $Y$ be a nonnegative continuous semimartingale and let $\sigma$ be a stopping time. Suppose that on $[0,\sigma)$,
    \begin{equation*}
        dY(t)=2\sqrt{a(t)Y(t)}\,dW(t)+b(t)\,dt, \qquad Y(0)=y_0>0,
    \end{equation*}
where $a,b$ are progressively measurable and locally integrable, and $a(t)>0$ for $t<\sigma$, almost surely. Define
    \begin{equation*}
        A(t)=\int_0^{t\wedge\sigma}a(s)\,ds,\qquad \theta(u)=\inf\{t\ge0:A(t)\ge u\},\qquad \widetilde Y(u)=Y(\theta(u)).
    \end{equation*}

\begin{lemma}[Random clock reduction]
    \label{lem:A-random-clock-reduction}
    Under the above assumptions, for $0\le u<A(\sigma-)$,
        \begin{equation*}
            d\widetilde Y(u) = 2\sqrt{\widetilde Y(u)}\,d\widetilde W(u) +\frac{b(\theta(u))}{a(\theta(u))}\,du
        \end{equation*}
    for a one-dimensional Brownian motion $\widetilde W$.
\end{lemma}

\begin{proof}
    Set
        \begin{equation*}
            N(t) = \int_0^{t\wedge\sigma}\sqrt{a(s)}\,dW(s).
        \end{equation*}
    Its quadratic variation is $A(t)$. Since $A$ is strictly increasing before $\sigma$,
        \begin{equation*}
            \widetilde W(u)=N(\theta(u)), \qquad 0\le u<A(\sigma-),
        \end{equation*}
    is Brownian motion in the time-changed filtration, up to that clock time. If necessary, extend it using an independent Brownian motion. The stochastic change-of-time formula gives
        \begin{equation*}
            \int_0^{\theta(u)} 2\sqrt{a(s)Y(s)}\,dW(s) = \int_0^u2\sqrt{\widetilde Y(v)}\,d\widetilde W(v).
        \end{equation*}
    The ordinary change of variables gives
        \begin{equation*}
            \int_0^{\theta(u)}b(s)\,ds = \int_0^u \frac{b(\theta(v))}{a(\theta(v))}\,dv.
        \end{equation*}
    These identities prove the assertion.
\end{proof}

\subsection{The product process and a squared-Bessel clock}

The product of the coordinates detects the first zero hit, also when several particles reach zero together. Positive collisions do not make this product vanish, so they do not require us to stop the argument. Define
    \begin{equation}
        \label{eq:A-boundary-notation}
        T_0:=\inf\{t\ge0:X_1(t)=0\}, \qquad T_-:=\inf\{t\ge0:X_1(t)<0\},
    \end{equation}
and
    \begin{equation*}
        E(t):=\prod_{i=1}^N X_i(t), \qquad G(t):=\sum_{i=1}^N\prod_{j\ne i}X_j(t),
    \end{equation*}
together with $d_A:=\kappa(\delta-N+1)$.

The following estimate will be used in the alternative proof below. It controls the integral which bounds the inverse product clock near zero.

The same factor $G$ appears in the drift and the quadratic variation of the product. Since $G>0$ before $T_0$, we can use it as the speed of the new clock. The resulting equation has the constant dimension parameter $d_A$.

For the product clock, write
    \begin{equation*}
        A(t)=\int_0^{t\wedge T_0}G(s)\,ds,\qquad \theta(u)=\inf\{t\ge0:A(t)>u\},\qquad u_*=\lim_{t\uparrow T_0}A(t).
    \end{equation*}
When $T_0=\infty$, the last limit is taken as $t\to\infty$. The inverse clock is used for $0\le u<u_*$.

\begin{lemma}[Product time change before the zero boundary]
    \label{lem:A-product-besq-time-change}
    Before $T_0$, the product satisfies
        \begin{equation}
            \label{eq:A-product-besq}
            dE(t)=2\sqrt{E(t)G(t)}\,dW(t)+d_A G(t)\,dt
        \end{equation}
    for a one-dimensional Brownian motion $W$. On an extension of the probability space if needed,
        \begin{equation*}
            E(t)=Z(A(t)),\qquad 0\le t<T_0,
        \end{equation*}
    where $Z$ is a squared Bessel diffusion of dimension $d_A$, started from $E(0)$ and stopped at its first zero.
\end{lemma}

\begin{proof}
    Before $T_0$, all coordinates are positive. Pairing the interaction terms in It\^o's formula gives
        \begin{equation*}
            \sum_{i<j}\frac{X_i+X_j}{X_i-X_j}\left(\frac1{X_i}-\frac1{X_j}\right)=-(N-1)\sum_i\frac1{X_i}.
        \end{equation*}
    The identity is used outside collision times, whose union has zero Lebesgue measure. Since $E\sum_iX_i^{-1}=G$,
        \begin{equation*}
            dE=2E\sum_i\frac{dB_i}{\sqrt{X_i}}+d_A G\,dt.
        \end{equation*}
    The martingale part has quadratic variation $4E^2\sum_iX_i^{-1}\,dt=4EG\,dt$, which proves \eqref{eq:A-product-besq}.

    The stopped process $Y(t)=E(t\wedge T_0)$ is nonnegative. Since $G>0$ before $T_0$, Lemma~\ref{lem:A-random-clock-reduction} with $\sigma=T_0$, $a=G$ and $b=d_AG$ gives
        \begin{equation*}
            dY(\theta(u)) =2\sqrt{Y(\theta(u))}\,d\widetilde W(u)+d_A\,du, \qquad 0\le u<u_*.
        \end{equation*}
    Extend the clock Brownian motion if necessary, and let $Z$ solve this squared Bessel equation from $E(0)$ up to its first zero. Pathwise uniqueness before zero gives $E(t)=Z(A(t))$ for $t<T_0$. For $d_A<0$, only the stopped diffusion is used.
\end{proof}
At this stage we do not assert that the available clock interval
reaches the first zero time of $Z$.

The time change alone does not show that zero is reached in finite time on the original time scale. We prove this by comparing the smallest particle with a squared Bessel process driven by the same Brownian motion.

\begin{lemma}[Finiteness of the physical hitting time]
    \label{lem:A-finite-zero-hit}
    If $d_A<2$, then $T_0<\infty$ almost surely.
\end{lemma}

\begin{proof}
    Before $T_0$, all particles are strictly positive. By ordering
    and non-stickiness, for Lebesgue-a.e. $t<T_0$,
    \begin{equation*}
        \frac{X_1(t)+X_j(t)}{X_1(t)-X_j(t)}
        \le -1, \qquad j=2,\ldots,N.
    \end{equation*}
    Consequently,
    \begin{equation*}
        dX_1(t)=2\sqrt{X_1(t)}\,dB_1(t)+b_1(t)\,dt,
    \end{equation*}
    where
    \begin{equation*}
        b_1(t)\le\kappa(\delta-N+1)=d_A.
    \end{equation*}

    Let $Z$ be a squared Bessel process of dimension $d_A<2$,
    driven by the same Brownian motion $B_1$, with
    $Z(0)=X_1(0)>0$, and stopped at its first zero time
    \begin{equation*}
        \zeta_0:=\inf\{t\ge0:Z(t)=0\}.
    \end{equation*}
    For $d_A<0$, the equation for $Z$ is understood only
    before $\zeta_0$.

    Put $D(t)=X_1(t)-Z(t)$. On $[0,T_0\wedge\zeta_0)$,
    \begin{equation*}
        d\langle D\rangle_t
        =4\left(\sqrt{X_1(t)}-\sqrt{Z(t)}\right)^2dt
        \le4|D(t)|\,dt.
    \end{equation*}
    As in Lemma~\ref{lem:A-zero-occupation}, the local time
    of $D$ at zero vanishes. Tanaka's formula therefore gives
    \begin{align*}
        D^+(t)
        ={}&2\int_0^t \mathbf 1_{\{D(s)>0\}}
        \left(\sqrt{X_1(s)}-\sqrt{Z(s)}\right)dB_1(s)\\
        &+\int_0^t \mathbf 1_{\{D(s)>0\}}
        (b_1(s)-d_A)\,ds,
    \end{align*}
    for $t<T_0\wedge\zeta_0$.

    Since $b_1(s)\le d_A$, the process $D^+$ is a
    nonnegative local supermartingale starting from zero.
    Hence $D^+(t)=0$ up to $T_0\wedge\zeta_0$, and therefore
    \begin{equation*}
        X_1(t)\le Z(t),
        \qquad t<T_0\wedge\zeta_0.
    \end{equation*}
    By continuity, $\zeta_0<T_0$ would imply
    \begin{equation*}
        0<X_1(\zeta_0)\le Z(\zeta_0)=0,
    \end{equation*}
    a contradiction. Thus $T_0\le\zeta_0$ almost surely.

    Since a squared Bessel process of dimension $d_A<2$
    hits zero in finite time almost surely, we conclude that
    $T_0<\infty$ almost surely.
\end{proof}

\subsection{Hitting zero for nonnegative dimension}
For $d_A\ge0$, Lemma~\ref{lem:A-indicator-nonnegative} shows that the process remains nonnegative. It remains to decide whether it reaches zero. The product representation excludes a zero hit when $d_A\ge2$, while Lemma~\ref{lem:A-finite-zero-hit} gives finite hitting time in the remaining range.

\begin{proposition}[Hitting zero without entrance into the negative half-line]
    \label{prop:A-hitting-zero-nonnegative-dimension}
    Assume $0<x_1<\cdots<x_N$. Then
        \begin{enumerate}[(i)]
            \item if $d_A\ge2$, then $T_0=\infty$ a.s. 
            \item if $0\le d_A<2$, then $T_0<\infty$, $T_-=\infty$ a.s.
        \end{enumerate}
\end{proposition}

\begin{proof}
    For $d_A\ge0$, Lemma~\ref{lem:A-indicator-nonnegative} gives $T_-=\infty$. Suppose first that $d_A\ge2$. Use the squared Bessel process $Z$ from Lemma~\ref{lem:A-product-besq-time-change}. If $T_0<\infty$, continuity and non-explosion imply
        \begin{equation*}
            A(T_0)=\int_0^{T_0}G(s)\,ds<\infty.
        \end{equation*}
    Since $E(t)=Z(A(t))$ for $t<T_0$, continuity gives
        \begin{equation*}
            Z(A(T_0))=E(T_0)=0,
        \end{equation*} 
    contradicting the polarity of zero for a squared Bessel process of dimension at least two. Thus $T_0=\infty$.

    If $0\le d_A<2$, Lemma~\ref{lem:A-finite-zero-hit} gives $T_0<\infty$, while nonnegative-chamber invariance gives $T_-=\infty$.
\end{proof}

\subsection{Directional exit from a zero block}

We now prove the new ingredient needed for the regime $d_A<0$.  Instead of studying the center of the whole zero block, we study its leftmost particle. The resulting drift estimate does not depend on the size of the zero block. We state the result for a general zero block at a stopping time, allowing negative particles on its left. At the first zero hit considered in this section, it will be applied with $\ell=0$.

\begin{lemma}[Immediate entrance into the negative half-line]
    \label{lem:A-zero-block-negative-exit}
    Let $X$ be the global non-sticky strong solution of \eqref{eq:A-particle-genuine}, and let $\tau<\infty$ be a stopping time such that
        \begin{equation*}
            X_1(\tau)\le\cdots\le X_\ell(\tau)<0 = X_{\ell+1}(\tau)=\cdots=X_{\ell+m}(\tau) < X_{\ell+m+1}(\tau)\le\cdots\le X_N(\tau),
        \end{equation*}
    for some
        \begin{equation*}
            0\le\ell\le N-m, \qquad m\ge1,
        \end{equation*}
    with the obvious conventions when $\ell=0$ or $\ell+m=N$. If
        \begin{equation}
            \label{eq:A-left-edge-negative-condition}
            \delta+2\ell+1-N<0,
        \end{equation}
    then
        \begin{equation}
            \label{eq:A-left-edge-immediate-negative}
            \inf\{t>0:X_{\ell+1}(\tau+t)<0\}=0 \qquad\text{a.s.}
        \end{equation}
\end{lemma}

\begin{proof}
    Shift time by $\tau$ and write
        \begin{equation*}
            Y(t):=X_{\ell+1}(\tau+t), \qquad \widetilde B_i(t):=B_i(\tau+t)-B_i(\tau).
        \end{equation*}
    Then $Y(0)=0$, and $\widetilde B$ is an $N$-dimensional Brownian motion with respect to the shifted filtration.  Set
        \begin{equation*}
            \sigma_-:=\inf\{t>0:Y(t)<0\}.
        \end{equation*}
    We prove that $\sigma_-=0$ almost surely.

    If $\ell\ge1$, define
        \begin{equation*}
            \rho_1 := \inf\left\{t>0: X_\ell(\tau+t)\ge\frac12X_\ell(\tau) \right\}.
        \end{equation*}
    Because $X_\ell(\tau)<0$, continuity gives $\rho_1>0$ almost surely.  If $\ell=0$, put $\rho_1=\infty$.  Also let
        \begin{equation*}
            \rho_2:=\inf\{t>0:Y(t)\ge1\}, \qquad \rho:=1\wedge\rho_1\wedge\rho_2.
        \end{equation*}
    Then $\rho>0$ almost surely. For $0\le t<\sigma_-\wedge\rho$, we have $0\le Y(t)<1$.  For every $j\le\ell$, the localization ensures $X_j(\tau+t)<0$, and hence
        \begin{equation*}
            \frac{|Y(t)|+|X_j(\tau+t)|}{Y(t)-X_j(\tau+t)} = 1.
        \end{equation*}
    For $j\ge\ell+2$, ordering gives
        \begin{equation*}
            X_j(\tau+t)\ge Y(t)\ge0.
        \end{equation*}
    Whenever $X_j(\tau+t)>Y(t)$,
        \begin{equation*}
            \frac{|Y(t)|+|X_j(\tau+t)|}{Y(t)-X_j(\tau+t)} = \frac{Y(t)+X_j(\tau+t)}{Y(t)-X_j(\tau+t)} \le-1.
        \end{equation*}
    The exceptional set on which $X_j=Y$ has zero Lebesgue measure by non-stickiness.  Therefore, $dt$-almost everywhere on $[0,\sigma_-\wedge\rho)$, the drift $b_Y$ of $Y$ satisfies
        \begin{align*}
            b_Y(t) &= \kappa\delta + \kappa\sum_{j\ne\ell+1} \frac{|Y(t)|+|X_j(\tau+t)|}{Y(t)-X_j(\tau+t)}\\
                   &\le \kappa\delta+\kappa\ell-\kappa(N-\ell-1)\\
                   &= \kappa(\delta+2\ell+1-N).
        \end{align*}
    By \eqref{eq:A-left-edge-negative-condition}, $c:=-\kappa(\delta+2\ell+1-N)>0$, so
        \begin{equation}
            \label{eq:A-left-edge-drift-bound}
            b_Y(t)\le-c
        \end{equation}
    for $dt$-almost every $t<\sigma_-\wedge\rho$. Before $\sigma_-$, $Y\ge0$, and hence
        \begin{equation*}
            dY(t)=2\sqrt{Y(t)}\,d\widetilde B_{\ell+1}(t)+b_Y(t)\,dt.
        \end{equation*}
    We now use the drift bound to rule out a positive waiting time before entrance into the negative half-line. After localization, such a waiting time would give a negative expectation for a nonnegative variable.
    Let $(\lambda_k)$ be a localizing sequence for the finite-variation part, chosen so that $\lambda_k>0$ almost surely and $\lambda_k\uparrow\infty$. Such a sequence exists because the singular drift integrals are finite on compact time intervals.  For $t>0$, set
        \begin{equation*}
            \theta_{t,k} := t\wedge\sigma_-\wedge\rho\wedge\lambda_k.
        \end{equation*}
    Since $0\le Y\le1$ before $\rho$, the stopped stochastic integral is a true martingale.  Taking expectations gives
        \begin{equation*}
            \mathbb E Y(\theta_{t,k}) = \mathbb E\int_0^{\theta_{t,k}}b_Y(s)\,ds.
        \end{equation*}
    By continuity, $Y(\theta_{t,k})\ge0$. Indeed, if the stopping occurs at $\sigma_-$, then $Y(\sigma_-)=0$; otherwise the process has not yet entered the negative half-line.  Using \eqref{eq:A-left-edge-drift-bound},
        \begin{equation*}
            0 \le \mathbb E Y(\theta_{t,k}) \le -c\,\mathbb E\theta_{t,k}.
        \end{equation*}
    Hence $\theta_{t,k}=0$ a.s.  Since $t>0$, $\rho>0$, and $\lambda_k>0$ almost surely, this forces $\sigma_-=0$ a.s.  which proves \eqref{eq:A-left-edge-immediate-negative}.
\end{proof}

\begin{rem}[Left edge versus block center]
    For a zero block of size $m$ with $\ell$ negative coordinates, Lemma~\ref{lem:zero-block-drift} gives the common drift assigned by the indicator equation at the face, $\kappa(\delta+2\ell+m-N)$. The preceding lemma instead uses the one-sided estimate $\kappa(\delta+2\ell+1-N)$ for the leftmost particle.  The latter does not depend on $m$.  This is the reason that no information about the multiplicity of the zero hit is needed to determine entrance into the negative half-line.
\end{rem}

\subsection{Negative dimension and the proof of the trichotomy}

It remains to consider the case $d_A<0$. We already know that zero is reached in finite time. We now show that the non-sticky solution cannot stay nonnegative on any interval immediately after the first zero. Several particles may reach zero together, so we have to allow a zero block of any size. Lemma~\ref{lem:A-zero-block-negative-exit} will be used at this point.

\begin{proposition}[Entrance into the negative half-line]
    \label{prop:A-negative-dimension-crossing}
    Assume $0<x_1<\cdots<x_N$ and $d_A=\kappa(\delta-N+1)<0$. Then $T_-=T_0<\infty$ almost surely. More precisely,
        \begin{equation*}
            \inf\{s>0:X_1(T_0+s)<0\}=0 \qquad\text{almost surely}.
        \end{equation*}
\end{proposition}

\begin{proof}
    Since $d_A<2$, Lemma~\ref{lem:A-finite-zero-hit} gives $T_0<\infty$ almost surely. By continuity and ordering, all coordinates are nonnegative at $T_0$, and the zero coordinates form an initial block. Let $M\in\{1,\ldots,N\}$ denote its size. On the event $\{M=m\}$, we have
        \begin{equation*}
            X_1(T_0)=\cdots=X_m(T_0)=0 < X_{m+1}(T_0)\le\cdots\le X_N(T_0),
        \end{equation*}
    where the positive part is omitted when $m=N$. The events $\{M=m\}$ belong to $\mathcal F_{T_0}$. On each of them, the argument of
    Lemma~\ref{lem:A-zero-block-negative-exit} applies at $\tau=T_0$ with $\ell=0$. The required condition is $\delta+1-N<0$, which follows from $d_A=\kappa(\delta-N+1)<0$ and $\kappa>0$. This condition does not depend on $m$. We therefore obtain
        \begin{equation*}
            \inf\{s>0:X_1(T_0+s)<0\}=0 \qquad\text{almost surely}.
        \end{equation*}
    Continuity gives $T_-\ge T_0$, since $X_1$ cannot become negative without first reaching zero. The last equality gives $T_-\le T_0$, because negative values occur arbitrarily soon after $T_0$. Thus $T_-=T_0<\infty$ almost surely.
\end{proof}

The equality $T_-=T_0$ concerns the infimum of the times when $X_1$ is negative. At the hitting time itself, continuity gives $X_1(T_0)=0$. The proposition says that every interval $(T_0,T_0+\varepsilon)$, $\varepsilon>0$, contains a time at which $X_1$ is negative. We can now complete the proof of Theorem~\ref{thm:A-complete-zero-trichotomy} stated in the Introduction.

\begin{proof}[Proof of Theorem~\ref{thm:A-complete-zero-trichotomy}]
    Proposition~\ref{prop:A-hitting-zero-nonnegative-dimension} gives parts~(i) and~(ii), which concern $d_A\ge0$. Proposition~\ref{prop:A-negative-dimension-crossing} gives part~(iii), including immediate entrance into the negative half-line when $d_A<0$.
\end{proof}

\subsection{Direct bounds for the first-zero survival probability}
\label{subsec:A-zero-tail}

We now estimate the probability that all coordinates remain strictly positive up to time $t$. This gives a quantitative complement to Theorem~\ref{thm:A-complete-zero-trichotomy}. Positive collisions are allowed before $T_0$. For $a>0$, $t>0$ and $z\ge0$, define
    \begin{equation}
        \label{eq:F}
        F_a(t,z)=\frac1{\Gamma(a)}\int_0^{z/(2t)}u^{a-1}e^{-u}\,du.
    \end{equation}
For $z>0$, direct differentiation gives
    \begin{equation}
        \label{eq:F-properties}
        \partial_tF_a=2z\partial_{zz}F_a+(2-2a)\partial_zF_a<0, \qquad \partial_zF_a>0.
    \end{equation}
Also $0\le F_a\le1$, $F_a(t,0)=0$, and $F_a(\varepsilon,z)\to1$ as $\varepsilon\downarrow0$ for $z>0$. For fixed $z>0$, the asymptotic below holds as $t\to\infty$; the upper bound holds for every $t>0$:
    \begin{equation}
        \label{eq:F-asymptotic}
        F_a(t,z)\sim\frac{z^a}{2^a\Gamma(a+1)}t^{-a},\qquad F_a(t,z)\le\min\left\{1,\frac{z^a}{2^a\Gamma(a+1)t^a}\right\}.
    \end{equation}
Below, $f(t)\asymp g(t)$ means two-sided positive constant bounds for all sufficiently large $t$, at a fixed starting point. The function $F_a$ gives the first-zero survival probability of a squared Bessel process of dimension $2-2a$. The next lemma gives an upper bound when the equation contains a random speed factor $b\ge1$ and an additional nonpositive drift.

\begin{lemma}[A scalar survival bound]
    \label{lem:scalar-comparison}
    Let $Z$ be continuous and positive before a stopping time $\theta$, with deterministic initial value $Z(0)=z>0$, and suppose $Z(t)\to0$ as $t\uparrow\theta$ whenever $\theta<\infty$. Suppose, before $\theta$,
        \begin{equation*}
            dZ=2\sqrt{bZ}\,dW+\bigl((2-2a)b+r\bigr)\,dt, \qquad a>0,\quad b\ge1,\quad r\le0,
        \end{equation*}
    with locally integrable coefficients and no explosion before $\theta$. Then $\mathbb P(\theta>t)\le F_a(t,Z(0))$. For a squared Bessel process of dimension $2-2a$, stopped at its first zero, this is an equality.
\end{lemma}

\begin{proof}
For fixed $t,\varepsilon>0$, the drift of $F_a(t-s+\varepsilon,Z(s))$ is $(b-1)\partial_tF_a+r\partial_zF_a\le0$. Localize before $\theta$ and use bounded convergence to obtain
    \begin{equation*}
        \mathbb E\bigl[\mathbf1_{\{\theta>t\}} F_a(\varepsilon,Z(t))\bigr] \le F_a(t+\varepsilon,Z(0)).
    \end{equation*}
The limiting value at a finite exit is zero. Let $\varepsilon\downarrow0$. When $b=1$ and $r=0$, the stopped process is a bounded martingale and the same argument gives equality.
\end{proof}

We can now prove Theorem~\ref{thm:A-positive-zero-survival}. For the upper bound, we apply the scalar estimate to $N$ times the geometric mean of the coordinates. For the lower bound, we use a bounded submartingale built from the product and the arithmetic mean.

\begin{proof}[Proof of Theorem \ref{thm:A-positive-zero-survival}]
    Use $E=\prod_iX_i$ and $G=e_{N-1}(X)$ from Lemma~\ref{lem:A-product-besq-time-change}. Before $T_0$, set
        \begin{equation*}
            R=NE^{1/N},\qquad a=\frac{G}{NE^{(N-1)/N}}\ge1.
        \end{equation*}
    The inequality is arithmetic--geometric mean applied to the summands of $G$. Equation~\eqref{eq:A-product-besq} and It\^o's formula give
        \begin{equation*}
            dR=2\sqrt{aR}\,dW+(2-2p_0)a\,dt.
        \end{equation*}
    Since $R$ vanishes at $T_0$, Lemma~\ref{lem:scalar-comparison} proves the upper bound.

    For the lower bound, let $h_0=E^{\nu_0}$ and $\overline X=N^{-1}\sum_iX_i$. Pairwise cancellation of interactions and $d_A+2\nu_0-2=0$ give
        \begin{equation*}
            dh_0=dM^0,\qquad d\overline X=dM^{\rm av}+\kappa\delta\,dt, \qquad d\langle M^{\rm av}\rangle_t=\frac4N\overline X\,dt, \qquad d\langle M^0,M^{\rm av}\rangle_t=4\nu_0h_0\,dt.
        \end{equation*}
    For fixed $t,\varepsilon>0$ and $0\le s\le t$, before $T_0$, set $Z_s=\overline X(s)+D_0(t-s+\varepsilon)$ and $Y_s=h_0(s)/Z_s^{p_0}$. The drift of $Y$ is
        \begin{equation*}
            Y_s\left[ \frac{p_0(D_0-\kappa\delta-4\nu_0)}{Z_s} +\frac{2p_0(p_0+1)\overline X(s)}{NZ_s^2}\right]ds =Y_s\frac{2p_0(p_0+1)\overline X(s)}{NZ_s^2}\,ds\ge0.
        \end{equation*}
    Arithmetic-geometric mean gives $0\le Y_s\le1$. Stop at
        \begin{equation*}
            \rho_m=\inf\{s\ge0:X_1(s)\le m^{-1} \text{ or }\textstyle\sum_iX_i(s)\ge m\},
        \end{equation*}
    which increases to $T_0$ and does not exclude positive collisions. The stopped process is a bounded submartingale. At a finite $T_0$ its numerator vanishes and its denominator is bounded away from zero by $D_0\varepsilon$. Hence
        \begin{equation*}
            \frac{h_0(0)}{[\overline{x}+D_0(t+\varepsilon)]^{p_0}} \le\mathbb E_x\!\left[\mathbf1_{\{T_0>t\}} \frac{h_0(t)}{[\overline X(t)+D_0\varepsilon]^{p_0}}\right] \le\mathbb P_x(T_0>t).
        \end{equation*}
    Let $\varepsilon\downarrow0$ and use $h_0(0)=g(x)^{p_0}$. All product and mean identities are valid through positive collisions: these are smooth symmetric functions, the singular drifts are locally integrable, and collision sets have zero occupation time. Finally, \eqref{eq:F-asymptotic} gives the matching powers.
\end{proof}

\begin{rem}
    If $d_A\ge2$, then $\mathbb P_x(T_0>t)=1$ by Theorem~\ref{thm:A-complete-zero-trichotomy}. If $d_A<0$, the bounds above also hold for $T_-$, since $T_-=T_0$ almost surely. At \(\kappa=1\) and \(N-1<\delta<N+1\), the first-zero problem is the first vanishing of a Wishart determinant. An exact survival formula in this range is given in \cite[Proposition~2.6]{MR2105711}, with the same power \(N(1-d_A/2)\).
\end{rem}

\section{Collisions}

Throughout this section, $X$ is the global non-sticky strong solution of \eqref{eq:A-particle-genuine}, with arbitrary $x\in\mathcal C_A$, $\kappa>0$, and $\delta\in\mathbb R$. We examine collisions at all positive times, including paths with earlier collisions and colliding initial configurations. Collisions at nonzero levels behave differently from those at zero. We first use the sum of squared pairwise distances to study the collapse of a block. Its equation can be compared locally with a squared Bessel equation. We use this comparison to study collisions away from zero, and then give an example of a triple collision at zero.

\subsection{Collision blocks and their squared distances}

\begin{dfn}[Maximal collision blocks]
    Let $I=\{\ell+1,\ldots,\ell+m\}$, where $m\ge2$. We call $I$ a maximal collision block of $z\in\mathcal C_A$ if
        \begin{equation*}
            z_{\ell+1}=\cdots=z_{\ell+m}=a,
        \end{equation*}
    and $z_\ell<a$ when $\ell\ge1$, while $a<z_{\ell+m+1}$ when $\ell+m<N$. The collision is called nonzero if $a\ne0$. No restriction is imposed on collisions among coordinates outside $I$.
\end{dfn}

For a contiguous block $I$ of size $m$, write
    \begin{equation}
        \label{eq:A-block-distance}
        \overline z_I=\frac1m\sum_{i\in I}z_i,\qquad
        S_I(z)=\sum_{\substack{i<j\\i,j\in I}}(z_j-z_i)^2
              =m\sum_{i\in I}(z_i-\overline z_I)^2.
    \end{equation}
The identity follows by expanding the squares. We write $S_I(t)$ for $S_I(X(t))$, and likewise use $\overline X_I(t)$. Thus $S_I=0$ corresponds to collapse of the entire block, whereas proper subblocks may collide when $S_I>0$.

We use deterministic neighborhoods of nonzero collision strata. For $0<\eta<L$ and $s\in\{-1,1\}$, define the relatively open subset of $\mathcal C_A$
    \begin{equation}
        \label{eq:A-block-chart}
        U_I^s(\eta,L)=\left\{z\in\mathcal C_A:
            \max_j|z_j|<L,\quad s z_i>\eta\ (i\in I),\quad
            \min_{\substack{i\in I\\j\notin I}}|z_i-z_j|>\eta
        \right\}.
    \end{equation}
The minimum over an empty set is $+\infty$. These neighborhoods allow both collisions among exterior particles and collisions of proper subblocks of $I$. We localize from a stopping time until the first subsequent exit from the chosen neighborhood. The separation from exterior particles allows us to control their contribution to the drift. For the interactions inside the block, the singular denominators cancel in the equation for $S_I$. This gives the local equation below.

\begin{lemma}[Local block equation]
    \label{lem:block_collision_local_structure}
    Fix $I$, $s$, $\eta$, and $L$ as above. On any stochastic interval in $U_I^s(\eta,L)$, the process $S_I$ has the decomposition
        \begin{equation}
            \label{eq:S-local-decomp}
            dS_I=dM_I+D_I\,dt,
            \qquad d\langle M_I\rangle_t=4S_I V_I\,dt,
        \end{equation}
    where, on $\{S_I>0\}$,
        \begin{equation}
            \label{eq:A-block-speed}
            V_I=4m\,
            \frac{\sum_{i\in I}|X_i|(X_i-\overline X_I)^2}
                 {\sum_{i\in I}(X_i-\overline X_I)^2},
        \end{equation}
    and
        \begin{align}
            D_I={}&(4+2\kappa m)
                    \sum_{\substack{i<j\\i,j\in I}}(|X_i|+|X_j|)
                    \label{eq:A-block-drift}\\
                 &+2\kappa m\sum_{i\in I}(X_i-\overline X_I)
                    \sum_{j\notin I}\frac{|X_i|+|X_j|}{X_i-X_j}.
                    \notag
        \end{align}
    The drift identities hold $dt$-almost everywhere, even if proper subblocks collide. Set
        \begin{equation}
            \label{eq:A-effective-dimension}
            d_m=(m-1)+\kappa\binom m2.
        \end{equation}
    There is a deterministic constant $C=C(N,m,\kappa,\eta,L)$ such that, whenever $0<S_I\le1$ in the neighborhood,
        \begin{equation}
            \label{eq:A-block-uniform-estimate}
            4m\eta\le V_I\le4mL,\qquad
            \left|\frac{D_I}{V_I}-d_m\right|\le C\sqrt{S_I}.
        \end{equation}
    Consequently, if the block converges to a common level $a\ne0$ while staying in such a neighborhood, then, along $\{S_I>0\}$,
        \begin{equation}
            \label{eq:DV-limits-block}
            V_I\longrightarrow4m|a|,\qquad
            D_I\longrightarrow4m|a|d_m,\qquad
            D_I/V_I\longrightarrow d_m.
        \end{equation}
\end{lemma}

\begin{proof}
    Set $Y_i=X_i-\overline X_I$ and $Q=\sum_{i\in I}Y_i^2$, so $S_I=mQ$. Apply It\^o's formula to the polynomial $m\sum_{i\in I}X_i^2-(\sum_{i\in I}X_i)^2$. The martingale part is
        \begin{equation*}
            dM_I=4m\sum_{i\in I}Y_i\sqrt{|X_i|}\,dB_i,
            \qquad d\langle M_I\rangle_t
                =16m^2\sum_{i\in I}|X_i|Y_i^2\,dt.
        \end{equation*}
    This yields \eqref{eq:A-block-speed}. The It\^o correction is $4(m-1)\sum_{i\in I}|X_i|$, while the constant drifts cancel because $\sum_iY_i=0$. At times with distinct coordinates, the internal interactions give
        \begin{equation*}
            2\kappa m\sum_{\substack{i<j\\i,j\in I}}
                \frac{(Y_i-Y_j)(|X_i|+|X_j|)}{X_i-X_j}
            =2\kappa m\sum_{\substack{i<j\\i,j\in I}}(|X_i|+|X_j|).
        \end{equation*}
    Non-stickiness makes the exceptional times Lebesgue-null, and local integrability of the genuine interaction terms justifies the integrated identity. Since $\sum_{i<j}(|X_i|+|X_j|)=(m-1)\sum_i|X_i|$, adding the exterior interactions gives \eqref{eq:A-block-drift}. Collisions of proper subblocks therefore introduce no additional terms.

    For fixed exterior coordinates, let
        \begin{equation*}
            F(u)=\sum_{j\notin I}\frac{|u|+|X_j|}{u-X_j}.
        \end{equation*}
    By contiguity of $I$ and ordering, every denominator has absolute value at least $\eta$ on the interval spanned by the block coordinates. Since these coordinates have the same sign, differentiation in $u$ shows that $F$ is Lipschitz there, with a deterministic constant depending only on $N,\eta,L$. Thus
        \begin{equation*}
            \left|\sum_{i\in I}Y_iF(X_i)\right|
            =\left|\sum_{i\in I}Y_i\bigl(F(X_i)-F(\overline X_I)\bigr)\right|
            \le C Q.
        \end{equation*}
    Since $\sum_i|X_i|=m|\overline X_I|$ in the neighborhood, we obtain
        \begin{equation*}
            D_I=4m|\overline X_I|d_m+O(S_I),\qquad
            |V_I-4m|\overline X_I||
                \le4m\max_{i\in I}|Y_i|\le4\sqrt{mS_I}.
        \end{equation*}
    The constants are uniform throughout the neighborhood. Formula \eqref{eq:A-block-speed} gives $4m\eta\le V_I\le4mL$, establishing \eqref{eq:A-block-uniform-estimate} and the stated limits. On $\{S_I=0\}$, choose any value of $V_I$ between these bounds; the quadratic-variation identity remains valid, and this set has zero occupation time. Where $S_I>0$, normalizing $M_I$ gives
        \begin{equation*}
            dS_I=2\sqrt{S_I V_I}\,dW+D_I\,dt
        \end{equation*}
    for a one-dimensional Brownian motion, with an independent extension if needed.
\end{proof}
After the time change with speed $V_I$, the squared distance has drift $D_I/V_I$. The estimate above allows us to keep this drift strictly above two near a nonzero collapse whenever $d_m>2$. We can then compare with a squared Bessel process which does not reach zero.

\begin{proposition}[Effective dimension and local nonattainment]
    \label{prop:effective_dimension_block_collision}
    Suppose $d_m>2$. For every neighborhood $U=U_I^s(\eta,L)$ there exists $r\in(0,1)$ with the following property. Let $\tau$ be a finite stopping time such that $X(\tau)\in U$ and $0<S_I(\tau)<r^2$, and define
        \begin{equation*}
            \sigma=\inf\{t\ge\tau:X(t)\notin U\ \text{or}\ S_I(t)\ge r^2\},\qquad
            T_I=\inf\{t\ge\tau:S_I(t)=0\}.
        \end{equation*}
    Then $\mathbb P(T_I\le\sigma,\ T_I<\infty)=0$. The same assertion holds restricted to the $\mathcal F_\tau$-measurable event on which these initial conditions hold. No assumption excludes earlier collisions or collisions of proper subblocks before $T_I$.
\end{proposition}

\begin{proof}
    Take $\varepsilon=(d_m-2)/2>0$. By \eqref{eq:A-block-uniform-estimate}, choose $r$ so that $D_I/V_I\ge2+\varepsilon$ whenever $0<S_I<r^2$ in $U$. Before $T_I\wedge\sigma$, introduce the clock
        \begin{equation*}
            A(t)=\int_\tau^t V_I(s)\,ds.
        \end{equation*}
    Lemma~\ref{lem:A-random-clock-reduction}, applied after shifting the starting time to $\tau$, gives
        \begin{equation*}
            d\widetilde S_I(u)=2\sqrt{\widetilde S_I(u)}\,d\widetilde W(u)
                +\widetilde b(u)\,du,\qquad \widetilde b(u)\ge2+\varepsilon.
        \end{equation*}
    With the same Brownian motion, couple this process to a squared Bessel process $Z$ of dimension $2+\varepsilon$, started from $Z(0)=S_I(\tau)>0$. One-dimensional comparison gives $\widetilde S_I\ge Z$ up to the stopped endpoint. This comparison holds first away from zero and then by localization, or directly by the standard square-root diffusion comparison theorem.

    Since $Z$ never hits zero, it is bounded away from zero on every compact clock interval. On $\{T_I\le\sigma,\ T_I<\infty\}$, the terminal clock is finite because $V_I\le4mL$ before that time. Continuity would then give $\widetilde S_I(A(T_I))=0$, contrary to comparison, even when $T_I=\sigma$. Localization by deterministic bounds on $\tau$ and $T_I$ completes the proof.
\end{proof}

\subsection{Multiple collisions away from zero}
The local result starts from a block with positive squared distance. We now use non-stickiness and continuity to find such a starting time before any possible nonzero collision at a positive time. This allows us to include initial configurations with coinciding particles.

\begin{theorem}[Nonzero block collisions at all positive times]
    \label{thm:no_positive_block_collision}
    Almost surely, no maximal nonzero collision block of size $m$ with $d_m>2$ occurs at any time $t>0$. In particular, for every $\kappa>0$,
        \begin{equation}
            \label{eq:A-no-nonzero-triples}
            \mathbb P\bigl(\exists\,t>0,\ i<j<k:
                X_i(t)=X_j(t)=X_k(t)\ne0\bigr)=0.
        \end{equation}
    If $\kappa>1$, nonzero two-particle collisions are also absent at all positive times. These assertions hold for every initial configuration in $\mathcal C_A$.
\end{theorem}

\begin{proof}
    Fix $I$ with $d_m>2$. Use the countable family of neighborhoods $U_I^s(1/n,L)$, where $n,L$ are positive integers with $1/n<L$, and choose one deterministic radius $r$ for each neighborhood as in Proposition~\ref{prop:effective_dimension_block_collision}. Apply that proposition from every positive rational time $q$, on the event
        \begin{equation*}
            X(q)\in U_I^s(1/n,L),\qquad 0<S_I(q)<r^2.
        \end{equation*}
    Outside a common null set, no such rational-time start is followed by a full block collapse before or at the corresponding exit.

    On the event of non-stickiness, suppose a maximal block $I$ collapses at a nonzero level at time $t>0$. Maximality and the nonzero level ensure that $X(t)$ lies in the relative interior of a neighborhood from the countable family. By continuity, for some $h\in(0,t)$,
        \begin{equation*}
            X(u)\in U_I^s(1/n,L),\qquad S_I(u)<r^2,
            \qquad t-h\le u\le t.
        \end{equation*}
    The set $\{u:S_I(u)=0\}$ lies in a pairwise collision set and hence has zero Lebesgue measure. There is therefore a time in $(t-h,t)$ with $S_I>0$, and by continuity such a rational time $q$ exists. Starting from $q$, the first full collapse occurs by time $t$ and before the exit, contradicting the local nonattainment assertion.

    There are finitely many contiguous blocks, so taking their union proves the first statement for all $t>0$ simultaneously. If $m\ge3$, then $d_m\ge2+3\kappa>2$; if $m=2$ and $\kappa>1$, then $d_2=1+\kappa>2$. Every collision is contained in a maximal block, proving the conclusions. The argument requires neither a first collision time for the full system nor a stopping time selected retrospectively from a collision event. Even when $S_I(0)=0$, non-stickiness supplies positive-distance starts arbitrarily close to zero, giving immediate departure and non-return for the excluded nonzero strata.
\end{proof}

\begin{rem}[Several collision blocks]
    The local calculation allows collisions among exterior particles. Thus the theorem also rules out nonzero multiple collisions simultaneous with disjoint collision blocks, whether at zero or at other nonzero levels. It does not rule out simultaneous collisions of two disjoint pairs when $0<\kappa<1$.
\end{rem}

\subsection{Two-particle collisions}

For a pair, the effective dimension is $d_2=1+\kappa$. The preceding argument excludes nonzero pair collisions when $\kappa>1$, while the critical case $\kappa=1$ requires a sharper drift estimate. The lower comparison below will exclude collisions for $\kappa\ge1$. For $0<\kappa<1$, we will use the upper comparison together with an estimate of the probability of leaving the neighborhood.

Fix $I=\{i,i+1\}$ and a neighborhood $U=U_I^s(\eta,L)$. Write
    \begin{equation*}
        G=X_{i+1}-X_i,\qquad S=G^2,\qquad K=|X_i|+|X_{i+1}|,\qquad d=1+\kappa.
    \end{equation*}
Let $\tau$ be a finite stopping time with $X(\tau)\in U$ and $S(\tau)>0$, and put
    \begin{equation*}
        \sigma=\inf\{t\ge\tau:X(t)\notin U\},\qquad T_I=\inf\{t\ge\tau:S(t)=0\}.
    \end{equation*}
The contribution of exterior particles to the gap drift is
    \begin{equation*}
        b_i=\kappa\sum_{j\notin I}\left(\frac{|X_{i+1}|+|X_j|}{X_{i+1}-X_j}-\frac{|X_i|+|X_j|}{X_i-X_j}\right).
    \end{equation*}
Lemma~\ref{lem:block_collision_local_structure}, applied with $m=2$, gives, before $T_I\wedge\sigma$,
    \begin{equation}
        \label{eq:A-squared-gap-equation}
        dS=2\sqrt{SV}\,dW+D\,dt,\qquad V=4K,\qquad D=4dK+2Gb_i.
    \end{equation}
We use the clock $A(t)=\int_\tau^tV(v)\,dv$, its inverse $\theta$, and $\widetilde S(u)=S(\theta(u))$. Lemma~\ref{lem:A-random-clock-reduction} gives
    \begin{equation}
        \label{eq:A-squared-gap-clock}
        d\widetilde S=2\sqrt{\widetilde S}\,d\widetilde W +\widetilde b\,du,\qquad \widetilde b=(D/V)\circ\theta
    \end{equation}
on the available clock interval. For a constant $C\ge0$, let $Z^-,Z^+$ start from $S(\tau)$ and use this same Brownian motion:
    \begin{equation}
        \label{eq:A-squared-gap-comparators}
        dZ^-=2\sqrt{Z^-}\,d\widetilde W+(d-CZ^-)\,du,\qquad dZ^+=2\sqrt{Z^+}\,d\widetilde W+d\,du.
    \end{equation}
The clock Brownian motion may be extended independently past its stopped endpoint.

\begin{lemma}[Squared-gap comparison]
    \label{lem:A-squared-gap-comparison}
    There is a deterministic $C\ge0$, depending only on $U$ and the parameters, such that before $T_I\wedge\sigma$,
        \begin{equation}
            \label{eq:A-squared-gap-drift-bounds}
            d-CS\le D/V\le d,\qquad 8\eta\le V\le8L.
        \end{equation}
    For this choice of $C$, the comparison $Z^-\le\widetilde S\le Z^+$ holds before the stopped clock endpoint. If $\kappa\ge1$, then
        \begin{equation*}
            \mathbb P(T_I\le\sigma,\ T_I<\infty)=0.
        \end{equation*}
    These assertions also hold on the $\mathcal F_\tau$-measurable event where the stated initial conditions hold.
\end{lemma}

\begin{proof}

    For each exterior coordinate $v=X_j$, the function $f_v(u)=(|u|+|v|)/(u-v)$ has derivative
        \begin{equation*}
            f_v'(u)=-\frac{\operatorname{sgn}(u)v+|v|}{(u-v)^2}\le0
        \end{equation*}
    on the interval $[X_i,X_{i+1}]$. The pair has a common nonzero sign in $U$, and adjacency places all exterior coordinates outside this interval. The neighborhood bounds also control $|f_v'|$ deterministically. Summation gives $-C_0G\le b_i\le0$. Since $K\ge2\eta$,
        \begin{equation*}
            \frac DV=d+\frac{Gb_i}{2K}
        \end{equation*}
    satisfies \eqref{eq:A-squared-gap-drift-bounds}, after adjusting the constant. Hence $d-C\widetilde S\le\widetilde b\le d$
in \eqref{eq:A-squared-gap-clock}.

    Square-root diffusion comparison, first away from zero and then by localization, gives $Z^-\le\widetilde S\le Z^+$. For the lower comparison, the drift bound used on $\{Z^->\widetilde S\}$ is
        \begin{equation*}
            d-CZ^-\le d-C\widetilde S\le\widetilde b.
        \end{equation*}
    If necessary, extend the clock Brownian motion independently beyond the stopped endpoint. A deterministic scaling and time change yield
        \begin{equation}
            \label{eq:A-damped-besq-representation}
            Z^-(u)=e^{-Cu}Y(h_C(u)),\qquad
            h_C(u)=
            \begin{cases}
                (e^{Cu}-1)/C,&C>0,\\
                u,&C=0,
            \end{cases}
        \end{equation}
    Here $Y$ is a squared Bessel process of dimension $d$, started from $S(\tau)$, with Brownian motion obtained by the corresponding deterministic time change. If $\kappa\ge1$, then $d\ge2$, so $Y$ never hits zero. Thus $Z^-$ remains strictly positive and is bounded away from zero on compact clock intervals. Any finite physical hit before or at $\sigma$ has finite clock time, since $V\le8L$, contradicting $\widetilde S\ge Z^-$. This also proves the assertion at the endpoint.
\end{proof}

\begin{theorem}[The critical threshold for nonzero pair collisions]
    \label{thm:two_particle_collisions}
    If $\kappa\ge1$, then, for every $\delta\in\mathbb R$ and every $x\in\mathcal C_A$,
        \begin{equation}
            \label{eq:A-no-nonzero-pairs}
            \mathbb P\bigl(\exists\,t>0,\ i<j:
                X_i(t)=X_j(t)\ne0\bigr)=0.
        \end{equation}
    In particular, an initially coinciding nonzero block separates immediately: almost surely its coordinates are pairwise distinct for all sufficiently small positive times.
\end{theorem}

\begin{proof}
    Theorem~\ref{thm:no_positive_block_collision} excludes nonzero collision blocks of size at least three. Any remaining nonzero collision must therefore involve a maximal adjacent pair. By Lemma~\ref{lem:A-squared-gap-comparison}, a positive gap inside $U_I^s(\eta,L)$ cannot reach zero before or at its neighborhood exit, at any finite time, for every $\kappa\ge1$.

    Apply the local assertion from every positive rational time in the countable family $U_I^s(1/n,L)$. If a maximal nonzero pair collided at time $t>0$, continuity would keep the path in one such neighborhood during a preceding interval. Non-stickiness and continuity then give a rational starting time in that interval with a positive gap, contradicting local nonattainment. Taking the finite union over adjacent pairs proves \eqref{eq:A-no-nonzero-pairs} simultaneously for all positive times, including after earlier zero collisions. If an initial block has nonzero common value, continuity keeps all its coordinates away from zero on a sufficiently short random interval; \eqref{eq:A-no-nonzero-pairs} then gives the final assertion.
\end{proof}

When $0<\kappa<1$, the upper comparator $Z^+$ has squared Bessel dimension $1+\kappa<2$. To establish accessibility before leaving the neighborhood, we first control dependence on the initial state. Write $\mathbb P_x$ for the law of the genuine solution started from $x$, and $\|z\|_1=\sum_j|z_j|$.

\begin{lemma}[Stability in the initial state]
    \label{lem:A-initial-state-stability}
    Let $X^x,X^y$ be the global genuine strong solutions started from $x,y\in\mathcal C_A$ and driven by the same Brownian motion. Then $R(t)=\|X^x(t)-X^y(t)\|_1$ is a nonnegative supermartingale. In particular, for every $T>0$ and $\varepsilon>0$,
        \begin{equation}
            \label{eq:A-initial-state-stability}
            \mathbb E R(T)\le\|x-y\|_1,\qquad
            \mathbb P\left(\sup_{0\le t\le T}R(t)\ge\varepsilon\right)
                \le\frac{\|x-y\|_1}{\varepsilon}.
        \end{equation}
    Thus, under this coupling, solutions converge uniformly on compact time intervals in probability as their initial states converge.
\end{lemma}

\begin{proof}
    The Tanaka argument of Lemma~\ref{lem:A-nonsticky-uniqueness} does not require identical initial states. The local times of the coordinate differences vanish, and the interaction contribution to $R$ is nonpositive. Hence
        \begin{equation*}
            R(t)=R(0)+M(t)-H(t),
        \end{equation*}
    where $M$ is a continuous local martingale started from zero and $H$ is increasing with $H(0)=0$. Local integrability follows from that of the genuine interaction terms. Thus $R$ is a nonnegative local supermartingale, hence a supermartingale. Optional stopping at $T\wedge\inf\{t\ge0:R(t)\ge\varepsilon\}$ gives the maximal inequality in \eqref{eq:A-initial-state-stability}.
\end{proof}
By stability and continuity, the probability of a very early exit from the neighborhood is small for initial states close to the collision configuration. The first zero time of the Bessel comparison scales as the square of the initial gap. Combining these estimates gives a collision before the exit with probability tending to one.

\begin{proposition}[Local subcritical accessibility]
    \label{prop:two_particle_subcritical_local_accessibility}
    Let $0<\kappa<1$ and $\delta\in\mathbb R$. Fix $x^0\in\mathcal C_A$ with a maximal nonzero collision block $I=\{i,i+1\}$, with common level $a\ne0$, and choose $U=U_I^s(\eta,L)$ containing $x^0$. For sufficiently small $g>0$, define $x^{(g)}\in U$ by
        \begin{equation*}
            x_i^{(g)}=a-g/2,\qquad x_{i+1}^{(g)}=a+g/2,\qquad
            x_j^{(g)}=x_j^0\quad(j\notin I).
        \end{equation*}
    Let
        \begin{equation*}
            \tau_I=\inf\{t\ge0:X_i(t)=X_{i+1}(t)\},\qquad
            \rho_U=\inf\{t\ge0:X(t)\notin U\}.
        \end{equation*}
    Then, for every $T>0$,
        \begin{equation}
            \label{eq:A-subcritical-local-accessibility}
            \lim_{g\downarrow0}
            \mathbb P_{x^{(g)}}(\tau_I\le T,\ \tau_I<\rho_U)=1.
        \end{equation}
    Thus the probability of an isolated nonzero collision before leaving $U$ tends to one as $g\downarrow0$ and is positive for all sufficiently small $g$. 
\end{proposition}

\begin{proof}
    Couple $X^{(g)}=X^{x^{(g)}}$ and $X^0=X^{x^0}$ using the same Brownian motion. Choose $\varepsilon>0$ such that the relative closed ball
        \begin{equation*}
            \{z\in\mathcal C_A:\|z-x^0\|_1\le2\varepsilon\}
            \subset U.
        \end{equation*}
    If $\sup_{v\le t}\|X^0(v)-x^0\|_1<\varepsilon$ and $\sup_{v\le t}\|X^{(g)}(v)-X^0(v)\|_1<\varepsilon$, then $X^{(g)}$ cannot exit $U$ by time $t$. Since $\|x^{(g)}-x^0\|_1=g$, Lemma~\ref{lem:A-initial-state-stability} gives
        \begin{equation*}
            \mathbb P_{x^{(g)}}(\rho_U\le t)
            \le p_0(t)+g/\varepsilon,\qquad
            p_0(t):=\mathbb P\left(\sup_{v\le t}\|X^0(v)-x^0\|_1\ge\varepsilon\right).
        \end{equation*}
    By continuity at zero, $p_0(t)\to0$ as $t\downarrow0$.

    For $X^{(g)}$, apply Lemma~\ref{lem:A-squared-gap-comparison} from time zero, stopped at $\tau_I\wedge\rho_U$. Extend its clock Brownian motion past the stopped endpoint and let $Z^{+,g}$ be the upper squared Bessel comparison, of dimension $1+\kappa$ and started at $g^2$. Write $\zeta_g$ for its first zero time. On $\{\tau_I>t,\ \rho_U>t\}$, the clock remains available through time $t$ and
        \begin{equation*}
            A(t)\ge8\eta t,\qquad
            0<\widetilde S(u)\le Z^{+,g}(u),\quad 0\le u\le A(t).
        \end{equation*}
    On this event, $\zeta_g>8\eta t$. Squared Bessel scaling gives $\zeta_g\overset{\mathrm{law}}=g^2\zeta_1$, where $\zeta_1$ denotes the first zero time for the same dimension started at one. Since $1+\kappa<2$, $\zeta_1<\infty$ almost surely. A union bound now yields, for every $t>0$,
        \begin{equation}
            \label{eq:A-subcritical-hit-bound}
            \mathbb P_{x^{(g)}}(\tau_I\le t,\ \tau_I<\rho_U)
            \ge1-p_0(t)-\frac g\varepsilon
                -\mathbb P\left(\zeta_1>\frac{8\eta t}{g^2}\right).
        \end{equation}
    This bound does not require independence of the comparison process and the neighborhood exit. For fixed $T>0$, use any $0<t\le T$, let $g\downarrow0$, and then let $t\downarrow0$. This proves \eqref{eq:A-subcritical-local-accessibility}. Before $\rho_U$, the pair has a common nonzero sign and remains separated from all exterior coordinates, so the collision is isolated and nonzero.
\end{proof}

\begin{rem}[Scope of local accessibility]
    The same proof applies to any family $x\to x^0$ in $\mathcal C_A$ with $x_i<x_{i+1}$: replace $g/\varepsilon$ by $\|x-x^0\|_1/\varepsilon$ and use the initial squared gap in the Bessel scaling. Thus sufficiently close initial configurations have positive probability of a collision before leaving the neighborhood. The result does not say that every initial configuration in the chamber can reach such a neighborhood. Collisions among exterior particles remain allowed.
\end{rem}

\subsection{Multiple collisions at zero}

Non-stickiness rules out positive collision occupation time, not instantaneous collisions. The preceding nonzero-collision argument also requires clock speed bounded away from zero. At a zero block collapse, \eqref{eq:A-block-speed} instead gives $V_I\to0$, so those estimates do not classify collisions at zero. The following example exhibits triple collisions at zero from a strictly ordered initial state.

\begin{proposition}[A sharp triple-zero threshold in a three-particle family]
    \label{prop:A-triple-zero-collision}
    Let $N=3$, $\delta=0$, $\kappa>0$, and $x=(-a,0,b)$ with $a,b>0$. The unique global genuine solution has the representation
        \begin{equation}
            \label{eq:A-three-particle-representation}
            X_1=-U,\qquad X_2=0,\qquad X_3=V,
        \end{equation}
    where $U,V$ are independent squared Bessel processes of dimension $2\kappa$, started at $a,b$. In particular, with
        \begin{equation*}
            \tau_{000}=\inf\{t>0:X_1(t)=X_2(t)=X_3(t)=0\},
        \end{equation*}
    one has
        \begin{equation*}
            \mathbb P(\tau_{000}<\infty)=
            \begin{cases}
                1,&0<\kappa<\tfrac12,\\
                0,&\kappa\ge\tfrac12.
            \end{cases}
        \end{equation*}
    When $0<\kappa<\tfrac12$, the triple collision occurs at a strictly positive finite time almost surely, although every pairwise collision set has zero Lebesgue measure.
\end{proposition}

\begin{proof}
    For the prescribed particle Brownian motion, take $W_1=-B_1$ and $W_2=B_3$, and let
        \begin{equation*}
            dU=2\sqrt U\,dW_1+2\kappa\,dt,\qquad
            dV=2\sqrt V\,dW_2+2\kappa\,dt.
        \end{equation*}
    The solutions are independent, nonnegative, and strong. Since $2\kappa>0$, both spend zero Lebesgue time at zero. Whenever $U,V>0$, the interactions of $-U$ with $0$ and $V$ are both $-1$, those of $V$ with $-U$ and $0$ are both $+1$, and those of the middle coordinate cancel. Thus the drifts in \eqref{eq:A-particle-genuine} are $-2\kappa,0,2\kappa$, respectively, and the signed Brownian motions give the required diffusion terms.

    The exceptional times lie in the Lebesgue-null union of the zero sets of $U,V$. Each defined interaction quotient has absolute value one, so the singular drift integrals are finite on compact intervals. Collisions with the pinned coordinate occur only in these zero sets; a collision of the exterior coordinates requires $U=V=0$. Thus \eqref{eq:A-three-particle-representation} is a non-sticky solution of the genuine equation and, by Theorem~\ref{thm:A-genuine-existence-uniqueness}, the unique solution.

    To determine when the three particles meet at zero, we consider the sum $Z=U+V$. Its drift is $4\kappa$ and its martingale part has quadratic variation $4Z\,dt$. Normalizing this martingale, with an independent Brownian extension on its zero set, gives
        \begin{equation*}
            dZ=2\sqrt Z\,dW+4\kappa\,dt,\qquad Z(0)=a+b>0.
        \end{equation*}
    Thus $Z$ is a squared Bessel process of dimension $4\kappa$. Started above zero, such a process hits zero almost surely in dimension below two and never hits zero in dimension at least two. Since $U,V\ge0$, the equality $Z=0$ is equivalent to $U=V=0$, which is exactly the triple collision. Its first hitting time is strictly positive by continuity and $Z(0)>0$.
\end{proof}

\begin{rem}
    The threshold in Proposition~\ref{prop:A-triple-zero-collision} applies to this three-particle family, not to all zero collision strata. Together with \eqref{eq:A-no-nonzero-triples}, it shows that degeneracy at zero allows collision behaviour excluded at nonzero levels.
\end{rem}

\section{The first collision time}
\label{sec:A-collision-time}

Let $X$ be the global non-sticky strong solution of \eqref{eq:A-particle-genuine}, started at
    \begin{equation*}
        x\in\mathcal C_A^\circ:=\{x\in\mathbb R^N:x_1<\cdots<x_N\}.
    \end{equation*}
The coordinates may have arbitrary signs; one may initially be zero. Define
    \begin{equation}
        \label{eq:A-collision-time}
        \tau_A=\inf\{t\ge0:X_i(t)=X_{i+1}(t) \text{ for some }1\le i<N\}.
    \end{equation}
Thus $\tau_A$ is the exit time from the open type-$A_{N-1}$ chamber. An individual coordinate hitting or crossing zero does not stop this clock. We use
    \begin{equation}
        \label{eq:A-collision-notation}
        \Delta(x)=\prod_{i<j}(x_j-x_i),\qquad S(x)=\sum_i|x_i|,\qquad M_N=\binom N2.
    \end{equation}
Our first goal in this section is to prove Theorem \ref{thm:A-collision-main}.

\subsection{Vandermonde identities in the signed chamber}
The Vandermonde product is positive before the first collision and vanishes at that collision. We use its logarithm to exclude collisions when $\kappa\ge1$, and the power $\Delta^{1-\kappa}$ to obtain the lower survival bound when $0<\kappa<1$.

For $x\in\mathcal C_A^\circ$, write $a_i=|x_i|$ and $g_i=\sum_{j\ne i}(x_i-x_j)^{-1}$. Set
    \begin{align}
        \mathcal A(x)&=\sum_{i<j}\frac{a_i+a_j}{(x_i-x_j)^2}, \label{eq:A-collision-A}\\
        \mathcal C(x)&=\sum_{i<j<k}\left[\frac{a_i}{(x_i-x_j)(x_i-x_k)}+\frac{a_j}{(x_j-x_i)(x_j-x_k)}+\frac{a_k}{(x_k-x_i)(x_k-x_j)}\right]. \label{eq:A-collision-C}
    \end{align}
An empty triple sum is zero. Each bracket is a second divided difference of the convex function $r\mapsto|r|$, so $\mathcal C\ge0$. A bracket vanishes whenever its three arguments belong to the same closed half-line. For $a<b<c$ with mixed signs the bracket is
\begin{equation*}
\begin{cases}
2(-a)/[(b-a)(c-a)],&a<0\le b<c,\\
2c/[(c-a)(c-b)],&a<b\le0<c.
\end{cases}
\end{equation*}
The two expressions agree at $b=0$. The term $\mathcal C$ comes from triples with mixed signs. For $0<\kappa<1$, choosing the power $1-\kappa$ cancels the term involving $\mathcal A$. The remaining drift is nonnegative, which will be used in the lower bound.

\begin{lemma}\label{lem:A-collision-vandermonde}
Let $\mathcal L$ be the generator in $\mathcal C_A^\circ$, and let
$M^{\log\Delta}$ be the local martingale part of $\log\Delta(X)$.
Before $\tau_A$,
\begin{align}
\mathcal L\log\Delta&=2(\kappa-1)\mathcal A+\kappa\mathcal C,
\label{eq:A-collision-log-drift}\\
d\langle M^{\log\Delta}\rangle_t
 &=4(\mathcal A+2\mathcal C)(X(t))\,dt.
\label{eq:A-collision-log-bracket}
\end{align}
For $0<\kappa<1$ and $h=\Delta^{1-\kappa}$,
\begin{equation}\label{eq:A-collision-h-drift}
\mathcal Lh=Vh,\qquad
V=(1-\kappa)(4-3\kappa)\mathcal C\ge0.
\end{equation}
These formulas also hold when one coordinate is zero and all gaps
are positive.
\end{lemma}
\begin{proof}
Grouping pairs and triples gives
\begin{equation*}
\sum_i a_i\sum_{j\ne i}\frac1{(x_i-x_j)^2}=\mathcal A,
\qquad \sum_i a_i g_i^2=\mathcal A+2\mathcal C,
\qquad
\sum_i g_i\sum_{j\ne i}\frac{a_i+a_j}{x_i-x_j}
=2\mathcal A+\mathcal C.
\end{equation*}
For the last identity, a triple contributes
$\sum_{r\in\{i,j,k\}}(a_r+a_i+a_j+a_k)/
\prod_{v\in\{i,j,k\}\setminus\{r\}}(x_r-x_v)$.
The constant numerator $a_i+a_j+a_k$ cancels, leaving the bracket in
\eqref{eq:A-collision-C}. Now use
$\partial_i\log\Delta=g_i$,
$\partial_{ii}\log\Delta=-\sum_{j\ne i}(x_i-x_j)^{-2}$ and
$\sum_i g_i=0$. For $q=1-\kappa$,
\begin{equation*}
\frac{\mathcal L\Delta^q}{\Delta^q}
=2q(q+\kappa-1)\mathcal A+q(\kappa+4q)\mathcal C
=q(4-3\kappa)\mathcal C.
\end{equation*}
Only smooth functions of the gaps are differentiated; no derivative
of $|x_i|$ is required.
\end{proof}

\subsection{Proof of the collision bounds}
For $\kappa\ge1$, the logarithmic argument below excludes collisions at zero as well as at nonzero levels. It uses identities which remain valid when one coordinate is zero and all gaps are positive.
\begin{proof}[Proof of Theorem~\ref{thm:A-collision-main}]
First suppose $\kappa\ge1$. By
\eqref{eq:A-collision-log-drift},
$\log\Delta(X)=\log\Delta(x)+M+A$ before $\tau_A$, where $A$ is
nondecreasing. At a finite collision, continuity and non-explosion
would force this process to be bounded above and tend to $-\infty$.
If $\langle M\rangle_{\tau_A-}<\infty$, the martingale has a finite
limit, which makes this impossible. If its quadratic variation is
infinite, its Brownian time change has infinite limsup, again
contradicting the upper bound for $M+A$. Hence no collision occurs.

Now assume $0<\kappa<1$. For the upper bound, let
$G_i=X_{i+1}-X_i$ and $A_i=|X_i|+|X_{i+1}|$.
Subtracting coordinate equations gives
\begin{equation}\label{eq:A-collision-gap}
dG_i=2\sqrt{A_i}\,d\beta_i+
       \left(\frac{2\kappa A_i}{G_i}+R_i(X)\right)dt,
\end{equation}
where
\begin{equation*}
R_i(x)=\kappa\sum_{j\notin\{i,i+1\}}
\left[\frac{|x_{i+1}|+|x_j|}{x_{i+1}-x_j}
     -\frac{|x_i|+|x_j|}{x_i-x_j}\right]\le0,
\qquad A_i\ge G_i.
\end{equation*}
Indeed, $f_z(y)=(|y|+|z|)/(y-z)$ is nonincreasing on every interval
avoiding $z$: away from $0,z$,
$f_z'(y)=-(z\operatorname{sgn}(y)+|z|)/(y-z)^2\le0$.
It is continuous at zero if $z\ne0$, and constant on each half-line
if $z=0$. Adjacency ensures that no other $x_j$ lies between the
endpoints of the gap.

Adjacent gaps have cross-variation
$d\langle G_i,G_{i+1}\rangle_t=-4|X_{i+1}|\,dt$;
nonadjacent gaps have zero cross-variation. For fixed
$t,\varepsilon>0$, put $r=t-s+\varepsilon$,
$F_i=F_q(r,G_i(s))$, and $\Phi_s=\prod_iF_i$.
Writing $F_{i,t}$ and $F_{i,z}$ for the derivatives, its drift is
\begin{align*}
&\sum_i\left[\left(\frac{A_i}{G_i}-1\right)F_{i,t}
+R_iF_{i,z}\right]\prod_{j\ne i}F_j-4\sum_{i=1}^{N-2}|X_{i+1}|F_{i,z}F_{i+1,z}
\prod_{j\notin\{i,i+1\}}F_j\le0.
\end{align*}
Thus this product is a bounded local supermartingale; independence
of the gaps is not needed. Localize at
\begin{equation}\label{eq:A-collision-localization}
\sigma_m=\inf\{s\ge0:\min_iG_i(s)\le m^{-1}
\text{ or }S(X(s))\ge m\}.
\end{equation}
These times increase to $\tau_A$. At a collision at least one
factor vanishes. Bounded convergence yields
\begin{equation*}
\mathbb E_x\!\left[\mathbf1_{\{\tau_A>t\}}
\prod_iF_q(\varepsilon,G_i(t))\right]
\le\prod_iF_q(t+\varepsilon,x_{i+1}-x_i).
\end{equation*}
Letting $\varepsilon\downarrow0$ proves the upper bound, and
\eqref{eq:F-asymptotic} also implies $\tau_A<\infty$ almost surely.

For the lower bound, apply Tanaka's formula to $S(X)=\sum_i|X_i|$.
The local times at zero vanish as in
Lemma~\ref{lem:A-zero-occupation}. With $\operatorname{sgn}(0)=0$,
\begin{align}
dS(X(t))&=dM^S(t)+b(X(t))\,dt,\qquad d\langle M^S\rangle_t=4S(X(t))\,dt,
\label{eq:A-collision-S}\\
b(x)&=\kappa\delta\sum_i\operatorname{sgn}(x_i)
+\kappa\sum_{i<j}|\operatorname{sgn}(x_i)-\operatorname{sgn}(x_j)|
\notag\le\kappa N|\delta|+2\kappa M_N.
\label{eq:A-collision-b}
\end{align}
The interaction term follows by pairing the drifts and using
\begin{equation*}
(\operatorname{sgn}(x_i)-\operatorname{sgn}(x_j))
\frac{|x_i|+|x_j|}{x_i-x_j}
=|\operatorname{sgn}(x_i)-\operatorname{sgn}(x_j)|.
\end{equation*}
For $h=\Delta^q$, homogeneity gives
\begin{equation}\label{eq:A-collision-cross-variation}
d\langle M^h,M^S\rangle_t
=4q h(X(t))\sum_iX_i(t)g_i(X(t))\,dt
=4p h(X(t))\,dt,
\end{equation}
because $\sum_i x_i g_i=M_N$. The homogeneity of $h$ suggests a denominator involving $S^p$. We add a positive time term so that the denominator stays away from zero, also when all particles meet at zero. The choice of $K$ will make the resulting bounded ratio a submartingale after localization. Set $Z_s=S(X(s))+K(t-s+\varepsilon)$ and $Y_s=h(X(s))/Z_s^p$. By \eqref{eq:A-collision-h-drift}, the drift of $Y$ is
\begin{equation}\label{eq:A-collision-submartingale-drift}
Y_s\left[V(X(s))+\frac{p(K-b(X(s))-4p)}{Z_s}
                    +\frac{2p(p+1)S(X(s))}{Z_s^2}\right]ds\ge0,
\end{equation}
since $K=\kappa N|\delta|+2\kappa M_N+4p$.
Each gap is at most $S$, so $0\le Y_s\le1$.
After stopping at $t\wedge\sigma_m$, this is a bounded
submartingale. Its limiting value at a collision is zero, since
$h=0$ and $Z_s\ge K\varepsilon$. Consequently,
\begin{equation*}
\frac{h(x)}{[S(x)+K(t+\varepsilon)]^p}
\le\mathbb E_x\!\left[\mathbf1_{\{\tau_A>t\}}
        \frac{h(X(t))}{[S(X(t))+K\varepsilon]^p}\right]
\le\mathbb P_x(\tau_A>t).
\end{equation*}
Let $\varepsilon\downarrow0$. This proves the lower bound, and
\eqref{eq:F-asymptotic} gives \eqref{eq:A-collision-powers}.
\end{proof}

\begin{rem}[Mixed signs and the two-particle case]
For mixed signs, $\mathcal C$ is generally nonzero; the term
$V=q(4-3\kappa)\mathcal C$ is retained in the lower-bound proof.
The estimates remain valid through individual zero crossings.
For $0<\kappa<1$ and $N=2$, they give
$\mathbb P_x(\tau_A>t)\asymp t^{-(1-\kappa)}$ for every $\delta$.
When also $\delta=0$ and $x_1\le0\le x_2$, with $x_1<x_2$,
the two magnitudes before collision are independent squared Bessel
processes of dimension $\kappa$. Their sum is a squared Bessel
process of dimension $2\kappa$, and the collision occurs at its first
zero. Hence $\mathbb P_x(\tau_A>t)=F_{1-\kappa}(t,x_2-x_1)$.
\end{rem}

\begin{rem}[Exit from the positive ordered chamber]
For $0<x_1<\cdots<x_N$, define the open positive ordered chamber
\begin{equation*}
\mathcal C_N^{+,\circ}
:=\{x\in\mathbb R^N:0<x_1<\cdots<x_N\}.
\end{equation*}
This is the interior of the closed chamber $\mathcal C_N^+$ used above.
Its exit time satisfies
\begin{equation*}
\tau_+=\inf\{t\ge0:X(t)\notin\mathcal C_N^{+,\circ}\}
       =T_0\wedge\tau_A.
\end{equation*}
If $d_A\ge2$, then $\tau_+=\tau_A$; if $\kappa\ge1$, then
$\tau_+=T_0$. If both inequalities hold, $\tau_+=\infty$ almost
surely. In the remaining regime the survival event is
$\{T_0>t,\tau_A>t\}$, and the two marginal survival probabilities
cannot in general be multiplied.
Nevertheless, if $0<\kappa<1$ and $d_A<2$, the event inclusion gives
    \begin{equation*}
        \mathbb P_x(\tau_+>t) \le\min\left\{F_{p_0}(t,Ng(x)),\, \prod_{i=1}^{N-1}F_q(t,x_{i+1}-x_i)\right\},
    \end{equation*}
where $p_0=N(1-d_A/2)$, $q=1-\kappa$, and $g(x)=(\prod_i x_i)^{1/N}$. This follows from Theorems~\ref{thm:A-positive-zero-survival} and~\ref{thm:A-collision-main}; it does not assert a sharp joint-survival exponent. Before $\tau_+$, the transformation $Y_i=\sqrt{X_i}$ gives a type-$B_N$ radial Dunkl equation with multiplicities $k_0=(d_A-1)/2$ and $k_1=\kappa/2$. For nonnegative multiplicities ($d_A\ge1$), boundary-hitting distributions are studied in \cite{MR2470522}. The first collision time $\tau_A$ is defined in the full ordered chamber, where individual coordinates may cross zero before a collision.
\end{rem}

\subsection{Collision survival for a nonnegative system}

Throughout this subsection, $0<\kappa<1$ and $q=1-\kappa$. We now consider an autonomous system of $n$ nonnegative particles. The term $\mathcal C$ then vanishes, and we can use the full Vandermonde product also for the upper bound. This gives upper and lower bounds with the same power of $t$. For $n\ge2$, put
    \begin{equation*}
        M_n=\binom n2,\quad p_n=M_nq,\quad \Delta_n(v)=\prod_{i<j}(v_j-v_i),\quad S_n(v)=\sum_i v_i, \quad R_n(v)=M_n\Delta_n(v)^{1/M_n}.
    \end{equation*}

Let $d\ge n-1$ and $0\le v_1<\cdots<v_n$. We consider the nonnegative solution of
    \begin{equation}
        \label{eq:block}
        dV_i=2\sqrt{V_i}\,dW_i + \kappa\left(d+\sum_{j\ne i}\frac{V_i+V_j}{V_i-V_j}\right)dt, \qquad V(0)=v.
    \end{equation}
Write
    \begin{equation*}
        \sigma_n=\inf\{t\ge0:V_i(t)=V_{i+1}(t) \text{ for some }1\le i<n\}, \qquad C_{n,d}=n\kappa d+4p_n.
    \end{equation*}
Individual visits to zero are allowed before $\sigma_n$.

\begin{lemma}[Collision survival in the nonnegative chamber]
    \label{lem:block}
    For every $t>0$,
        \begin{equation}
            \label{eq:block-collision-bounds}
            \frac{\Delta_n(v)^q}{[S_n(v)+C_{n,d}t]^{p_n}} \le\mathbb P_v(\sigma_n>t) \le F_{p_n}(t,R_n(v)).
        \end{equation}
    In particular, $\mathbb P_v(\sigma_n>t)\asymp t^{-p_n}$.
\end{lemma}

\begin{proof}
Existence and uniqueness follow from
Theorem~\ref{thm:A-genuine-existence-uniqueness}, with $N=n$ and
$\delta=d$. Lemma~\ref{lem:A-indicator-nonnegative} gives
nonnegativity because $d\ge n-1$. Only the stopped process is used.
Write $h=\Delta_n^q$, $S=S_n$ and
\begin{equation*}
A=\sum_{i<j}\frac{V_i+V_j}{(V_j-V_i)^2}.
\end{equation*}
All coordinates are nonnegative, so $\mathcal C=0$ in the
Vandermonde identities. Consequently
\begin{equation*}
\mathcal Lh=0,\qquad
d\langle h\rangle_t=4q^2h^2A\,dt,\qquad
d\langle h,S\rangle_t=4p_nh\,dt,
\end{equation*}
while $\mathcal LS=n\kappa d$ and
$d\langle S\rangle_t=4S\,dt$.
As in \eqref{eq:A-collision-submartingale-drift},
\begin{equation*}
\frac{h(V(s))}{[S(V(s))+C_{n,d}(t-s+\varepsilon)]^{p_n}}
\end{equation*}
is bounded by one and has nonnegative drift; it vanishes at the first
collision. This proves the lower bound by localization and bounded
convergence.

For the upper bound, we use the geometric mean of all pairwise gaps. The factor $M_n$ in $R_n$ makes the speed coefficient $b$ at least one, as required by Lemma~\ref{lem:scalar-comparison}. It\^o's formula for $R=R_n(V)$ gives
    \begin{equation}
        \label{eq:block-R}
        dR=2\sqrt{bR}\,d\beta+(2-2p_n)b\,dt, \qquad b=\frac{Rq^2A}{p_n^2}\ge1.
    \end{equation}
Indeed,
    \begin{equation*}
        A\ge\sum_{i<j}\frac1{V_j-V_i} \ge M_n\Delta_n(V)^{-1/M_n},
    \end{equation*}
by arithmetic-geometric mean. Apply Lemma~\ref{lem:scalar-comparison}. The functions $h$ and $R$ are smooth wherever all gaps are positive, including when $V_1=0$; no stopping at a single zero visit is needed.
\end{proof}

\bibliographystyle{abbrv}

\end{document}